# Dynamics of a Predator-Prey Model with Allee Effect and Interspecific Competition

Peng Lina*          Xie Jianhang

## Abstract

This paper primarily discusses the dynamical properties of a class of Lotka-Volterra models featuring the Allee effect and interspecific competition within the predator population. The constructed models employ Holling II and Holling I response functions for the predator, respectively.The existence of boundary equilibrium points under various parameter conditions and internal equilibrium points under specific parameter conditions is discussed. The equilibrium points of the system may be stable or unstable nodes, saddle points, saddle-nodes, or cusp points with a codimension of 2. The parameter conditions under which internal equilibrium points possess one zero eigenvalue and two non-zero eigenvalues, one zero eigenvalue and a pair of purely imaginary eigenvalues, or two zero eigenvalues and one non-zero eigenvalue are analyzed.

# 1 Introduction

Direct relationships exist between populations, such as predation, competition, cooperation, and parasitism. Due to the complexity and diversity of populations and their ecological environments, indirect factors such as interspecific cooperation, the fear effect, and the Allee effect should also be considered. This paper primarily focuses on the Allee effect and interspecific competition.

The Allee effect describes the relationship between population density and individual fitness, taking into account the effects of both intraspecific cooperation and intraspecific competition on population dynamics. In the 20th century, Stephens[1] summarized the Allee effect, classifying it into the following two types:

Weak Allee effect: this refers to a situation where the population growth rate increases as population density increases, resulting in a positive growth rate. A typical expression for the weak Allee effect is:

$$f(x) = rx\left(1 - \frac{x}{K}\right)\frac{x}{x+A}, \tag{1}$$

where $\frac{x}{x+A}$ represents the Allee term.

Strong Allee effect: this refers to a situation where the population growth rate increases as population density increases; however, if the population density falls below the Allee effect threshold, the growth rate becomes negative. A common expression for the strong Allee effect is:

$$f(x) = rx\left(1 - \frac{x}{K}\right)(x - m), \tag{2}$$

where $m$ represents the Allee threshold.

Sasmal[2] proposed a Lotka-Volterra predator-prey model in which prey growth rates are reduced due to the fear effect and prey species are influenced by the Allee effect:

$$\begin{cases} \dfrac{du}{dt} = ru\left(1 - \dfrac{u}{K}\right)(u - \theta)\dfrac{1}{1 + fv} - auv \\ \dfrac{dv}{dt} = a\alpha uv - mv \end{cases}, \tag{3}$$

and investigated its dynamical behavior. By transforming the original system into a dimensionless system through variable transformations, the study found that the original system possesses three boundary equilibrium points as well as a positive equilibrium point that appears under certain conditions. The types of equilibrium points and their stability were analyzed, revealing that the origin is always locally asymptotically stable.

Predatory relationships are common in nature, while competitive relationships are further divided into interspecific and intraspecific competition. Interspecific competition refers to the competition between populations for shared resources such as water, sunlight, and soil nutrients; the outcome of such competition is often asymmetric, meaning that only one party gains an advantage among the competing parties. Intraspecific competition can regulate the size of a biological population; for example, if a beetle lays too many eggs, it will consume some of them. However, this is constrained by population density; that is, the higher the population density, the more intense the intraspecific competition. Furthermore, for some animals that live in family groups, intraspecific competition can also mitigate the harmful effects of long-term inbreeding on the population, helping it adapt to the environment.

Ikbal[3] proposed a Lotka-Volterra three-species predator-prey model involving interspecific and intraspecific competition between two predator populations:

$$\begin{cases} \dfrac{dp}{dt} = rp - \dfrac{rp^2}{K} - \dfrac{\alpha_1 P H_1}{1+\theta P} - \alpha_2 P H_2 \\ \dfrac{dH_1}{dt} = \dfrac{e_1 \alpha_1 P H_1}{1+\theta P} - g_1 H_1^2 - \beta_1 H_1 H_2 - d_1 H_1 \\ \dfrac{dH_2}{dt} = e_2 \alpha_2 P H_2 - g_2 H_2^2 - \beta_2 H_1 H_2 - d_2 H_2 \end{cases} \quad (4)$$

and analyze the existence of positive feedback, the existence of equilibrium points, and the stability of the system near equilibrium points.

Based on the model developed by the aforementioned researchers, the following three-species predator-prey model, which incorporates interspecific competition among predator populations, is established on the basis of the Lotka-Volterra model with the Allee effect:

$$\begin{cases} \dfrac{dx}{dt} = rx\left(1 - \dfrac{x}{K}\right)(x-m) - \dfrac{c_1 x y_1}{1+\theta x} - c_2 x y_2 \\ \dfrac{dy_1}{dt} = -s_1 y_1 + \dfrac{\alpha_1 x y_1}{1+\theta x} - \beta_1 y_1 y_2 \\ \dfrac{dy_2}{dt} = -s_2 y_2 + \alpha_2 x y_2 - \beta_2 y_1 y_2 \end{cases} \quad (5)$$

The initial conditions are

$$x(0) > 0 \ ; y_1(0) \geq 0 \ ; y_2(0) \geq 0 \ . \tag{6}$$

where $r_1$ is the intrinsic growth rate of the prey, $K$ is the environmental carrying capacity, $m$ is the Allee effect threshold for the prey, and $c_i$ represents the predation rate of the predator population $i$, $\theta$ is the semi-saturation constant, $s_i$ is the natural mortality rate of the predator population $i$, $\alpha_i$ is the predator birth rate generated by the predator population $i$ during predation, and $\beta_i$ is the disturbance rate experienced by the predator population $i$ during interspecific competition.

All parameters in Model (5) are positive. Using the following transformation:

$$\tilde{x} = \frac{x}{K}, \tilde{y}_1 = \frac{c_1 y_1}{rK}, \tilde{y}_2 = \frac{c_2 y_2}{rK}, \tau = rKt, \tilde{m} = \frac{m}{K}, \tilde{\theta} = \theta K,$$

$$\tilde{s}_1 = \frac{s_1}{rK}, \tilde{s}_2 = \frac{s_2}{rK}, \tilde{\alpha}_1 = \frac{\alpha_1}{r}, \tilde{\alpha}_2 = \frac{\alpha_2}{r}, \tilde{\beta}_1 = \frac{\beta_1}{c_2}, \tilde{\beta}_2 = \frac{\beta_2}{c_1}.$$

Replace $x, y_1, y_2, t, m, \theta, s_1, s_2, \alpha_1, \alpha_2, \beta_1, \beta_2$ with $\tilde{x}, \tilde{y}_1, \tilde{y}_2, \tau, \tilde{m}, \tilde{\theta}, \tilde{s}_1, \tilde{s}_2, \tilde{\alpha}_1, \tilde{\alpha}_2, \tilde{\beta}_1, \tilde{\beta}_2$ to obtain the following model

$$\begin{cases} \dfrac{dx}{dt} = x(1-x)(x-m) - \dfrac{xy_1}{1+\theta x} - xy_2 \\ \dfrac{dy_1}{dt} = -s_1 y_1 + \dfrac{\alpha_1 x y_1}{1+\theta x} - \beta_1 y_1 y_2 \\ \dfrac{dy_2}{dt} = -s_2 y_2 + \alpha_2 x y_2 - \beta_2 y_1 y_2 \end{cases}, \tag{7}$$

where the threshold range for the strong Allee effect is $0 < m < 1$.

## 2 Dynamical Properties of Equilibrium Points

In predator-prey models, equilibrium points hold significant ecological importance. Mathematically, an equilibrium point refers to a point where the rates of change of all variables in the model's system of differential equations are zero; that is, at this point, the population densities of both the predator and the prey reach equilibrium and no longer change. Therefore, the equilibrium points in the predator-prey model not only hold mathematical analytical value but also provide an important scientific foundation for the sustainable management of ecosystems and the prediction of population dynamics.

Based on the biological significance of the variables in the system, find

$$\Omega = \{(x, y_1, y_2) \mid x(0) > 0, y_1(0) \geq 0, y_2(0) \geq 0\},$$

we will solve for the equilibrium point of (7), that is, find the solution to the following system of equations

$$\begin{cases} x(1-x)(x-m) - \dfrac{xy_1}{1+\theta x} - xy_2 = 0 \\ -s_1 y_1 + \dfrac{\alpha_1 x y_1}{1+\theta x} - \beta_1 y_1 y_2 = 0 \\ -s_2 y_2 + \alpha_2 x y_2 - \beta_2 y_1 y_2 = 0 \end{cases} \quad (8)$$

**Theorem 1.** The boundary equilibrium points of the above system (7) are as follows:

(i) The equilibrium point of the free boundary of the preyed-upon population, $E_1 = (1, 0, 0)$, and $E_2 = (m, 0, 0)$.

(ii) The equilibrium point on the boundary where Predator Population 1 is free

$$E_3 = (\dfrac{s_2}{\alpha_2}, 0, \dfrac{(\alpha_2 - s_2)(s_2 - m\alpha_2)}{\alpha_2^2}).$$

(iii) The equilibrium point on the boundary where Predator Population 2 is free

$$E_4 = (\dfrac{s_1}{\alpha_1 - \theta s_1}, \dfrac{\alpha_1(\alpha_1 - \theta s_1 - s_1)(s_1 - m\alpha_1 + m\theta s_1)}{(\alpha_1 - \theta s_1)^2}, 0).$$

Proof. For the system of equations (8), consider the case where the density of at least one of the two predator populations is zero.

When $y_1 = y_2 = 0$, we have

$$x(1-x)(x-m) = 0,$$

which yields $x_1 = 1$, and $x_2 = m$. Thus, the system has boundary equilibrium points $E_1 = (1, 0, 0)$ and $E_2 = (m, 0, 0)$.

When $y_1 = 0$, we have

$$\begin{cases} x(1-x)(x-m)-xy_2=0 \\ -s_2y_2+\alpha_2xy_2=0 \end{cases}, \tag{9}$$

Since $y_2 \neq 0$, we have $-s_2+\alpha_2 x=0$. Solving yields $x=\dfrac{s_2}{\alpha_2}$. Substituting $x$ into the first term of (9)

gives $y_2=\dfrac{(\alpha_2-s_2)(s_2-m\alpha_2)}{\alpha_2^2}$. Thus, the system has boundary equilibrium points

$$E_3=(\dfrac{s_2}{\alpha_2},0,\dfrac{(\alpha_2-s_2)(s_2-m\alpha_2)}{\alpha_2^2}).$$

When $y_2=0$, we have

$$\begin{cases} x(1-x)(x-m)-\dfrac{xy_1}{1+\theta x}=0 \\ -s_1y_1+\dfrac{\alpha_1 xy_1}{1+\theta x}=0 \end{cases}, \tag{10}$$

Since $y_1 \neq 0$, we have $-s_1+\dfrac{\alpha_1 x}{1+\theta x}=0$, which yields $x=\dfrac{s_1}{\alpha_1-\theta s_1}$,

Substituting $x$ into the first term of (10) yields

$$y_1=\dfrac{\alpha_1(\alpha_1-\theta s_1-s_1)(s_1-m\alpha_1+m\theta s_1)}{(\alpha_1-\theta s_1)^2}.$$

Thus, the system has a boundary equilibrium point

$$E_4=(\dfrac{s_1}{\alpha_1-\theta s_1},\dfrac{\alpha_1(\alpha_1-\theta s_1-s_1)(s_1-m\alpha_1+m\theta s_1)}{(\alpha_1-\theta s_1)^2},0). \qquad \square$$

**Theorem 2.** If $\beta_1=\beta_2, \theta=1, s_1=s_2, m=\dfrac{2s_1}{\beta_1}$, the internal equilibrium points of the system (7)

above are as follows:

(i) When $0<s_1<\dfrac{\beta_1}{2}$ and $\alpha_2 > \dfrac{s_1^2+s_1\beta_1+1}{2\beta_1}$, the system has no internal equilibrium points.

(ii) When $0<s_1<\dfrac{\beta_1}{2}$ and $\alpha_2 = \dfrac{s_1^2+s_1\beta_1+1}{2\beta_1}$, the system has a double equilibrium point

$E_5=(x_5, y_1^{*5}, y_2^{*5})$.

(iii) When $0 < s_1 < \dfrac{\beta_1}{2}$ and $0 < \alpha_2 < \dfrac{s_1^2 + s_1\beta_1 + 1}{2\beta_1}$, the system has two equilibrium points:

$E_6 = (x_6, y_1^{*6}, y_2^{*6})$ and $E_7 = (x_7, y_1^{*7}, y_2^{*7})$.

Proof. Solving for the internal equilibrium points of system (7) is equivalent to finding the exact solution to the system of equations (8). We obtain:

$$y_1 = \dfrac{-s_2 + \alpha_2 x}{\beta_2}, y_2 = \dfrac{-s_1 - s_1\theta x + \alpha_2 x}{\beta_1(1+\theta x)}.$$

where $x$ is a positive root of the following univariate cubic equation (10)

$$\beta_1\beta_2\theta x^3 - \beta_1\beta_2(\theta + \theta m - 1)x^2 + (\theta m - 1 - m + \alpha_2\beta_1 + \alpha_2\beta_2 - \beta_2 s_1\theta)x + \beta_1\beta_2 m - \beta_1 s_2 - \beta_2 s_1 = 0.$$

Given the computational complexity of finding roots of a cubic equation, we will now provide a sufficient condition for the existence of a positive root of the above equation.

Let $\beta_1 = \beta_2, \theta = 1, s_1 = s_2, m = \dfrac{2s_1}{\beta_1}$. Since $0 < m < 1$, then $0 < s_1 < \dfrac{\beta_1}{2}$.

Under the above conditions, Equation (10) can be written as

$$\beta_1^2 x^2 - 2s_1\beta_1 x + 2\alpha_2\beta_1 - \beta_1 s_1 - 1 = 0. \tag{11}$$

Next, we examine the case of positive roots of Equation (11).

Let the discriminant of equation (11) be $\Delta = 4s_1^2\beta_1^2 - 4\beta_1^2(2\alpha_2\beta_1 - s_1\beta_1 - 1)$.

Let $\Delta < 0$, then $s_1^2 - 2\alpha_2\beta_1 + s_1\beta_1 + 1 < 0$, so $\alpha_2 > \dfrac{s_1^2 + s_1\beta_1 + 1}{2\beta_1}$, since $\beta_1 > 0$. Considering the sign of $f(s_1) = s_1^2 + \beta_1 s_1 + 1$, let $\Delta_{s_1} = \beta_1^2 - 4s_1^2$. Since $0 < s_1 < \dfrac{\beta_1}{2}$, then $\Delta_{s_1} > 0$. Meanwhile, $f(0) = 1$, and the axis of symmetry of $f(s_1)$ is $\hat{x} = -\dfrac{\beta_1}{2} < 0$. Therefore, when $s_1 > 0$, $f(s_1)$ is always greater than 0.

Based on the sign of $f(s_1)$, when $0 < s_1 < \dfrac{\beta_1}{2}$ and $\alpha_2 > \dfrac{s_1^2 + s_1\beta_1 + 1}{2\beta_1}$, we have $\Delta < 0$, and the system has no internal equilibrium points at this point.

Let $\Delta = 0$. Then, when $0 < s_1 < \dfrac{\beta_1}{2}$ and $\alpha_2 = \dfrac{s_1^2 + s_1\beta_1 + 1}{2\beta_1}$, equation (11) has a unique positive root

$x_5$, and $x_5 = \dfrac{s_1}{\beta_1}$. Therefore, system (7) has a unique internal equilibrium point

$E_5 = (x_5, y_1^{*5}, y_2^{*5})$, where $y_1^{*5} = \dfrac{s_1(\alpha_2 - \beta_1)}{\beta_1^2}$, $y_2^{*5} = \dfrac{s_1(\alpha_2 - s_1 - \beta_1)}{\beta_1(1+\theta x)}$.

Let $\Delta > 0$. Then, when $0 < s_1 < \dfrac{\beta_1}{2}$ and $0 < \alpha_2 < \dfrac{s_1^2 + s_1\beta_1 + 1}{2\beta_1}$, equation (11) has two positive roots

$x_6 = \dfrac{s_1 - \sqrt{s_1^2 - 2\alpha_2\beta_1 - \beta_1 s_1 - 1}}{\beta_1}$ and $x_7 = \dfrac{s_1 + \sqrt{s_1^2 - 2\alpha_2\beta_1 - \beta_1 s_1 - 1}}{\beta_1}$, and $x_6 < x_7$. That is, system

(7) has two internal equilibrium points $E_6 = (x_6, y_1^{*6}, y_2^{*6})$ and $E_7 = (x_7, y_1^{*7}, y_2^{*7})$, where

$y_1^{*i} = \dfrac{-s_1 + \alpha_2 x_i}{\beta_1}$, $y_2^{*i} = \dfrac{-s_1 - s_1\theta x_i + \alpha_2 x_i}{\beta_1(1+x)}$ and $i = 6, 7$. $\square$

The stability of an equilibrium point reflects an ecosystem's ability to return to equilibrium after being disturbed. A stable equilibrium point indicates that the system possesses self-regulatory capacity and can return to its initial state following external disturbances, whereas an unstable equilibrium point may foreshadow population extinction. Below, we will analyze the type and stability of equilibrium points using the eigenvalues of the system's Jacobian matrix at the equilibrium points.

If $E = (x, y_1, y_2)$ is an equilibrium point of system (7), then the Jacobian matrix of the system at this equilibrium point is

$$J(E) = \begin{pmatrix} -3x^2 + 2(m+1)x - m - \dfrac{y_1}{(1+\theta x)^2} - y_2 & -\dfrac{x}{1+\theta x} & -x \\ \dfrac{\alpha_1 y_1}{(1+\theta x)^2} & -s_1 + \dfrac{\alpha_1 x}{1+\theta x} - \beta_1 y_2 & -\beta_1 y_1 \\ \alpha_2 y_2 & -\beta_2 y_2 & -s_2 + \alpha_2 x - \beta_2 y_1 \end{pmatrix}.$$

By analyzing the eigenvalues of the Jacobian matrix at the equilibrium point, we obtain the following theorem.

**Theorem 3.** For the boundary equilibrium point $E_1 = (1, 0, 0)$,

(i) When $\dfrac{\alpha_1}{1+\theta} < s_1$ and $\alpha_2 < s_2$, the equilibrium point $E_1$ is a locally stable point.

(ii) When $\frac{\alpha_1}{1+\theta} > s_1$, $\alpha_2 \neq s_2$, or $\frac{\alpha_1}{1+\theta} \neq s_1$, $\alpha_2 > s_2$, the equilibrium point $E_1$ is a saddle point.

Proof. The Jacobian matrix of system (7) at the equilibrium point $E_1 = (1,0,0)$ is

$$J(E_1) = \begin{pmatrix} m-1 & -\frac{1}{1+\theta} & -1 \\ 0 & \frac{\alpha_1}{1+\theta} - s_1 & 0 \\ 0 & 0 & \alpha_2 - s_2 \end{pmatrix}.$$

$J(E_1)$ with eigenvalues $\lambda_1 = m-1$, $\lambda_2 = \frac{\alpha_1}{1+\theta} - s_1$, and $\lambda_3 = \alpha_2 - s_2$. For the strong Allee effect, we have $0 < m < 1$, so $\lambda_1 < 0$. Next, we consider the signs of $\lambda_2$ and $\lambda_3$. When $\frac{\alpha_1}{1+\theta} < s_1$ and $\alpha_2 < s_2$, we have $\lambda_2 < 0$ and $\lambda_3 < 0$, so the boundary equilibrium point $E_1$ is a locally stable node. When $\frac{\alpha_1}{1+\theta} > s_1$, $\alpha_2 \neq s_2$, or $\frac{\alpha_1}{1+\theta} \neq s_1$, and $\alpha_2 > s_2$, we have $\lambda_2 > 0$, $\lambda_3 \neq 0$, or $\lambda_2 \neq 0$, and $\lambda_3 > 0$; the boundary equilibrium point $E_1$ is a saddle point. □

**Theorem 4.** For the boundary equilibrium point $E_2 = (m,0,0)$

(i) When $\frac{\alpha_1 m}{1+\theta m} > s_1$ and $\alpha_2 m > s_2$, the equilibrium point $E_2$ is an unstable point.

(ii) When $\frac{\alpha_1 m}{1+\theta m} < s_1$, $\alpha_2 m \neq s_2$, or $\frac{\alpha_1 m}{1+\theta m} \neq s_1$, and $\alpha_2 m = s_2$, the equilibrium point $E_2$ is a saddle point.

Proof. The Jacobian matrix of system (7) at the equilibrium point $E_2 = (m,0,0)$ is:

$$J(E_2) = \begin{pmatrix} -m^2 + m & -\frac{m}{1+\theta m} & -m \\ 0 & \frac{\alpha_1 m}{1+\theta m} - s_1 & 0 \\ 0 & 0 & \alpha_2 m - s_2 \end{pmatrix}.$$

$J(E_1)$ with eigenvalues $\lambda_1 = -m^2 + m$, $\lambda_2 = \frac{\alpha_1 m}{1+\theta m} - s_1$, and $\lambda_3 = \alpha_2 m - s_2$. For the strong Allee effect, we have $0 < m < 1$, hence $\lambda_1 > 0$. Below, we consider the signs of $\lambda_2$ and $\lambda_3$. When $\frac{\alpha_1 m}{1+\theta m} > s_1$ and $\alpha_2 m > s_2$, we have $\lambda_2 > 0$ and $\lambda_3 > 0$; thus, the boundary equilibrium point $E_2$ is an unstable node. When $\frac{\alpha_1 m}{1+\theta m} < s_1, \alpha_2 m \neq s_2$, or $\frac{\alpha_1 m}{1+\theta m} \neq s_1, \alpha_2 m = s_2$, we have $\lambda_2 < 0, \lambda_3 \neq 0$, or $\lambda_2 \neq 0$, $\lambda_3 < 0$; thus, the boundary equilibrium point $E_2$ is a saddle point.

□

**Theorem 5.** For the boundary equilibrium point $E_3 = (\frac{s_2}{\alpha_2}, 0, \frac{(\alpha_2 - s_2)(s_2 - m\alpha_2)}{\alpha_2^2})$

(i) When $s_1 > \frac{\alpha_1 s_2}{\alpha_2 + \theta s_2} - \frac{\beta_1(\alpha_2 - s_2)(s_2 - m\alpha_2)}{\alpha_2^2}$ and $s_2 > \frac{(m+1)\alpha_2}{4}$, the equilibrium point $E_3$ is a locally asymptotically stable node.

(ii) When $s_1 > \frac{\alpha_1 s_2}{\alpha_2 + \theta s_2} - \frac{\beta_1(\alpha_2 - s_2)(s_2 - m\alpha_2)}{\alpha_2^2}$ and $0 < s_2 < \frac{(m+1)\alpha_2}{4}$, the equilibrium point $E_3$ is a saddle point.

(iii) When $0 < s_1 < \frac{\alpha_1 s_2}{\alpha_2 + \theta s_2} - \frac{\beta_1(\alpha_2 - s_2)(s_2 - m\alpha_2)}{\alpha_2^2}$ and $s_2 > \frac{(m+1)\alpha_2}{4}$, the equilibrium point $E_3$ is a saddle point.

(iv) When $0 < s_1 < \frac{\alpha_1 s_2}{\alpha_2 + \theta s_2} - \frac{\beta_1(\alpha_2 - s_2)(s_2 - m\alpha_2)}{\alpha_2^2}$ and $0 < s_2 < \frac{(m+1)\alpha_2}{4}$, the equilibrium point $E_3$ is a locally asymptotically stable node.

(v) When $s_1 = \frac{\alpha_1 s_2}{\alpha_2 + \theta s_2} - \frac{\beta_1(\alpha_2 - s_2)(s_2 - m\alpha_2)}{\alpha_2^2}$ and $s_2 \neq \frac{(m+1)\alpha_2}{4}$, the equilibrium point $E_3$ is a saddle-node point of codimension 1.

Proof.  The Jacobian matrix of system (7) at the equilibrium point

$$E_3 = (x_3, y_1^{*3}, y_2^{*3}) = (\frac{s_2}{\alpha_2}, 0, \frac{(\alpha_2 - s_2)(s_2 - m\alpha_2)}{\alpha_2^2})$$

is

$$J(E_3) = \begin{pmatrix} \dfrac{s_2(m\alpha_2 + \alpha_2 - 4s_2)}{\alpha_2^2} & -\dfrac{s_2}{\alpha_2 + \theta s_2} & -\dfrac{s_2}{\alpha_2} \\ 0 & -s_1 + \dfrac{\alpha_1 s_2}{\alpha_2 + \theta s_2} - \dfrac{\beta_1(\alpha_2 - s_2)(s_2 - m\alpha_2)}{\alpha_2^2} & 0 \\ \dfrac{(\alpha_2 - s_2)(s_2 - m\alpha_2)}{\alpha_2} & -\dfrac{\beta_2(\alpha_2 - s_2)(s_2 - m\alpha_2)}{\alpha_2^2} & 0 \end{pmatrix}.$$

The characteristic equation is

$|\lambda I - J(E_3)|$
$= (\lambda + s_1 - \dfrac{\alpha_1 s_2}{\alpha_2 + \theta s_2} + \dfrac{\beta_1(\alpha_2 - s_2)(s_2 - m\alpha_2)}{\alpha_2^2})(\lambda^2 - \dfrac{s_2(m\alpha_2 + \alpha_2 - 4s_2)}{\alpha_2^2}\lambda + \dfrac{s_2(\alpha_2 - s_2)(s_2 - m\alpha_2)}{\alpha_2^2}).$

The eigenvalues of $J(E_3)$ are $\lambda_1 = -s_1 + \dfrac{\alpha_1 s_2}{\alpha_2 + \theta s_2} - \dfrac{\beta_1(\alpha_2 - s_2)(s_2 - m\alpha_2)}{\alpha_2^2}$, $\lambda_2$, and $\lambda_3$, which are

roots of the following equation

$$\lambda^2 - \dfrac{s_2(m\alpha_2 + \alpha_2 - 4s_2)}{\alpha_2^2}\lambda + \dfrac{s_2(\alpha_2 - s_2)(s_2 - m\alpha_2)}{\alpha_2^2}.$$

Since $y_2^{*3} = \dfrac{(\alpha_2 - s_2)(s_2 - m\alpha_2)}{\alpha_2^2} > 0$, it follows that

$$\lambda_2 \cdot \lambda_3 = \dfrac{s_2(\alpha_2 - s_2)(s_2 - m\alpha_2)}{\alpha_2^2} > 0.$$

(i) Since $s_1 > \dfrac{\alpha_1 s_2}{\alpha_2 + \theta s_2} - \dfrac{\beta_1(\alpha_2 - s_2)(s_2 - m\alpha_2)}{\alpha_2^2}$, then

$$\lambda_1 = -s_1 + \dfrac{\alpha_1 s_2}{\alpha_2 + \theta s_2} - \dfrac{\beta_1(\alpha_2 - s_2)(s_2 - m\alpha_2)}{\alpha_2^2} < 0,$$

and since $s_2 > \dfrac{(m+1)\alpha_2}{4}$, then: $\lambda_2 + \lambda_3 = \dfrac{s_2(m\alpha_2 + \alpha_2 - 4s_2)}{\alpha_2^2} < 0$  Thus, the eigenvalues of $J(E_4)$ are $\lambda_1 < 0$, $\lambda_2 < 0$, and $\lambda_3 < 0$.

Therefore, when $s_1 > \dfrac{\alpha_1 s_2}{\alpha_2 + \theta s_2} - \dfrac{\beta_1(\alpha_2 - s_2)(s_2 - m\alpha_2)}{\alpha_2^2}$ and $s_2 > \dfrac{(m+1)\alpha_2}{4}$, the equilibrium point

$E_3$ is locally asymptotically stable.

(ii) Since $s_1 > \dfrac{\alpha_1 s_2}{\alpha_2 + \theta s_2} - \dfrac{\beta_1(\alpha_2 - s_2)(s_2 - m\alpha_2)}{\alpha_2^2}$, then

$$\lambda_1 = -s_1 + \frac{\alpha_1 s_2}{\alpha_2 + \theta s_2} - \frac{\beta_1(\alpha_2 - s_2)(s_2 - m\alpha_2)}{\alpha_2^2} < 0,$$

furthermore, since $s_2 < \frac{(m+1)\alpha_2}{4}$, we have $\lambda_1 + \lambda_2 = \frac{s_2(m\alpha_2 + \alpha_2 - 4s_2)}{\alpha_2^2} > 0$.

Thus, the eigenvalues of $J(E_4)$ are $\lambda_1 < 0$, $\lambda_2 > 0$, and $\lambda_3 > 0$. Therefore, when $s_1 > \frac{\alpha_1 s_2}{\alpha_2 + \theta s_2} - \frac{\beta_1(\alpha_2 - s_2)(s_2 - m\alpha_2)}{\alpha_2^2}$ and $s_2 < \frac{(m+1)\alpha_2}{4}$, the equilibrium point $E_3$ is a saddle point.

(iii) Since $0 < s_1 < \frac{\alpha_1 s_2}{\alpha_2 + \theta s_2} - \frac{\beta_1(\alpha_2 - s_2)(s_2 - m\alpha_2)}{\alpha_2^2}$, then

$$\lambda_1 = -s_1 + \frac{\alpha_1 s_2}{\alpha_2 + \theta s_2} - \frac{\beta_1(\alpha_2 - s_2)(s_2 - m\alpha_2)}{\alpha_2^2} > 0,$$

and since $s_2 > \frac{(m+1)\alpha_2}{4}$, then: $\lambda_1 + \lambda_2 = \frac{s_2(m\alpha_2 + \alpha_2 - 4s_2)}{\alpha_2^2} < 0$.

Thus, the eigenvalues of $J(E_4)$ are $\lambda_1 > 0$, $\lambda_2 < 0$, and $\lambda_3 < 0$. Therefore, when $0 < s_1 < \frac{\alpha_1 s_2}{\alpha_2 + \theta s_2} - \frac{\beta_1(\alpha_2 - s_2)(s_2 - m\alpha_2)}{\alpha_2^2}$ and $s_2 > \frac{(m+1)\alpha_2}{4}$, the equilibrium point $E_3$ is a saddle point.

(iv) Since $0 < s_1 < \frac{\alpha_1 s_2}{\alpha_2 + \theta s_2} - \frac{\beta_1(\alpha_2 - s_2)(s_2 - m\alpha_2)}{\alpha_2^2}$, then

$$\lambda_1 = -s_1 + \frac{\alpha_1 s_2}{\alpha_2 + \theta s_2} - \frac{\beta_1(\alpha_2 - s_2)(s_2 - m\alpha_2)}{\alpha_2^2} > 0,$$

and since $0 < s_2 < \frac{(m+1)\alpha_2}{4}$, then: $\lambda_2 + \lambda_3 = \frac{s_2(m\alpha_2 + \alpha_2 - 4s_2)}{\alpha_2^2} > 0$.

Thus, the eigenvalues of $J(E_4)$ are $\lambda_1 > 0$, $\lambda_2 > 0$, and $\lambda_3 > 0$. Therefore, when $0 < s_1 < \frac{\alpha_1 s_2}{\alpha_2 + \theta s_2} - \frac{\beta_1(\alpha_2 - s_2)(s_2 - m\alpha_2)}{\alpha_2^2}$ and $0 < s_2 < \frac{(m+1)\alpha_2}{4}$, the equilibrium point $E_3$ is an unstable node.

(v) When

$$s_1 = \frac{\alpha_1 s_2}{\alpha_2 + \theta s_2} - \frac{\beta_1(\alpha_2 - s_2)(s_2 - m\alpha_2)}{\alpha_2^2}, \quad \alpha_1 > \frac{\beta_1(\alpha_2 - s_2)(s_2 - m\alpha_2)(\alpha_2 + \theta s_2)}{\alpha_2^2 s_2} > 0$$

we have $\lambda_1 = 0$.

Furthermore, since $s_2 \neq \frac{(m+1)\alpha_2}{4}$, we have $\lambda_2 + \lambda_3 = \frac{s_2(m\alpha_2 + \alpha_2 - 4s_2)}{\alpha_2^2} \neq 0$. Thus, the eigenvalues of $J(E_4)$ are $\lambda_1 = 0$, $\lambda_2 \neq 0$, and $\lambda_3 \neq 0$.

Perform the transformation

$$\begin{cases} u = x - x_3 \\ v = y_1 \\ z = y_2 - y_2^{*3} \end{cases},$$

by translating the equilibrium point to the origin and performing a Taylor expansion at the origin, we obtain

$$\begin{cases} \dot{u} = a_{100}u + a_{010}v + a_{001}z + a_{200}u^2 + a_{110}uv + a_{101}uz + a_{300}u^3 + a_{210}u^2v + o(\|(u,v,z)\|^4) \\ \dot{v} = b_{010}v + b_{110}uv + b_{011}vz + b_{210}u^2v + o(\|(u,v,z)\|^4) \\ \dot{z} = c_{100}u + c_{010}v + c_{001}z + c_{101}uz + c_{011}vz + o(\|(u,v,z)\|^4) \end{cases}, \quad (12)$$

where $a_{100} = -3x_3^2 + 2(m+1)x_3 - m - y_2^{*3}, a_{010} = -\frac{x_3}{1+\theta x_3} a_{010} = -\frac{x_3}{1+\theta x_3}$, $a_{001} = -x_3$,

$a_{200} = m + 1 - 3x_3$, $a_{110} = -\frac{1}{(1+\theta x_3)^2}$, $a_{101} = -1$, $a_{300} = -1$, $a_{210} = \frac{1}{(1+\theta x_3)^3}$,

$b_{001} = -s_1 + \frac{\alpha_1 x_3}{1+\theta x_3} - \beta_1 y_2^{*3}$, $b_{110} = \frac{\alpha_1}{(1+\theta x_3)^2}$, $b_{011} = -\beta_1$, $b_{210} = -\frac{\alpha_1}{(1+\theta x_3)^3}$,

$c_{100} = \alpha_2 y_2^{*3}$, $c_{010} = -\beta_2 y_2^{*3}$, $c_{001} = -s_2 + \alpha_2 x_3$, $c_{101} = \alpha_2$, $c_{011} = -\beta_2$.

To transform the linearization matrix of system (12) into standard form, the following transformation is applied to system (12)

$$\begin{pmatrix} u \\ v \\ z \end{pmatrix} = \begin{pmatrix} w_{11} & w_{12} & w_{13} \\ w_{21} & w_{22} & w_{23} \\ w_{31} & w_{32} & w_{33} \end{pmatrix} \begin{pmatrix} u_1 \\ v_1 \\ z_1 \end{pmatrix},$$

Including: $w_{11} = 1$, $w_{12} = 1$, $w_{13} = 1$, $w_{21} = 0$, $w_{22} = 0$,

$$w_{23} = \alpha_1 - \frac{s_1(\alpha_2 + \theta s_2)}{s_2} - \frac{\beta_1(\alpha_2 - s_2)(s_2 - m\alpha_2)(\alpha_2 + \theta s_2) + s_2(m\alpha_2 + \alpha_2 - 4s_2)(\alpha_2 + \theta s_2)}{\alpha_2^2 s_2}$$

$$w_{31} = \frac{m\alpha_2 + \alpha_2 - 4s_2}{\alpha_2} \; , \; w_{32} = 0 \; , \; w_{33} = 0.$$

Result

$$\begin{cases} \dot{u}_1 = d_{200}u_1^2 + d_{110}u_1v_1 + d_{101}u_1z_1 + o(\|(u_1, v_1, z_1)\|^3) \\ \dot{v}_1 = \lambda_2 v_1 + e_{200}u_1^2 + e_{110}u_1v_1 + e_{101}u_1z_1 + e_{020}v_1^2 + e_{011}v_1z_1 + e_{002}z_1^2 + o(\|(u_1, v_1, z_1)\|^3) \; , \\ \dot{z}_1 = \lambda_3 z_1 + f_{101}u_1z_1 + f_{011}v_1z_1 + f_{002}z_1^2 + o(\|(u_1, v_1, z_1)\|^3) \end{cases} \quad (13)$$

where $d_{200} = c_{101}$, $d_{110} = c_{101}$, $d_{101} = c_{101} + c_{011}w_{23}$, $e_{200} = a_{200} + a_{101}w_{21} - c_{101}$,

$e_{110} = 2a_{200} + a_{101}w_{21} - c_{101}$, $e_{101} = 2a_{200} + a_{110}w_{23} + a_{101}w_{21} - b_{110} - b_{011}w_{21} - c_{101} - c_{011}w_{23}$

$e_{020} = a_{200}$, $e_{011} = 2a_{200} + a_{110}w_{23} - b_{110}$, $e_{002} = a_{200} + a_{110}w_{23} - b_{110}$, $f_{101} = b_{110} + b_{011}w_{21}$,

$f_{011} = b_{110}$, $f_{002} = b_{110}$

Let the center manifold of (13) be: $\begin{pmatrix} v_1 \\ z_1 \end{pmatrix} = \begin{pmatrix} h_1(u_1) \\ h_2(u_1) \end{pmatrix} = \begin{pmatrix} \sigma_1 u_1^2 + o(|u_1|^3) \\ \sigma_2 u_1^2 + o(|u_1|^3) \end{pmatrix}$, then

$$\begin{pmatrix} 2\sigma_1 u_1 + \cdots \\ 2\sigma_2 u_1 + \cdots \end{pmatrix} (d_{200}u_1^2 + d_{110}u_1 \cdot \sigma_1 u_1^2 + d_{101}u_1 \cdot \sigma_2 u_1^2 + \cdots)$$
$$= \begin{pmatrix} \lambda_2 \cdot \sigma_1 u_1^2 + e_{200}u_1^2 + e_{110}u_1 \cdot \sigma_1 u_1^2 + e_{101}u_1 \cdot \sigma_2 u_1^2 + e_{020}(\sigma_1 u_1^2)^2 + e_{011} \cdot \sigma_1 u_1^2 \cdot \sigma_2 u_1^2 + e_{002}(\sigma_2 u_1^2)^2 + \cdots \\ \lambda_3 \cdot \sigma_2 u_1^2 + f_{101}u_1 \cdot \sigma_2 u_1^2 + f_{011} \cdot \sigma_1 u_1 \cdot \sigma_2 u_1^2 + f_{002}(\sigma_2 u_1^2)^2 + \cdots \end{pmatrix}.$$

By comparing the coefficients, we obtain: $\lambda_2 \sigma_1 + e_{200} = 0$, $\lambda_3 \sigma_2 = 0$, so

$$\sigma_1 = -\frac{e_{200}}{\lambda_2} \; , \; \hat{\sigma}_2 = 0,$$

hence the central manifold is

$$\begin{pmatrix} v_1 \\ z_1 \end{pmatrix} = \begin{pmatrix} h_1(u_1) \\ h_2(u_1) \end{pmatrix} = \begin{pmatrix} -\frac{e_{200}}{\lambda_2}u_1^2 + o(|u_1|^3) \\ o(|u_1|^3) \end{pmatrix}.$$

The system (13) is equivalent near the origin to the system:

$$\dot{u}_1 = g_{200}u_1^2 + g_{110}u_1^3 + o(\|(u_1, v_1, z_1)\|^3) \; ,$$

where $g_{200} = d_{200} = c_{101} = \alpha_2 > 0$, $g_{110} = -\frac{e_{200}}{\lambda_2} d_{110}$. By the theorem, we know that $E_3$ is a saddle-node point of codimension 1. □

**Theorem 6.** For boundary Equilibrium Points $E_4$

(i) When $s_2 > \frac{\alpha_2 s_1 - \alpha_1 \beta_2(\alpha_1 - \theta s_1 - s_1)(s_1 - m\alpha_1 + m\theta s_1)}{\alpha_1 - \theta s_1}$ and $\alpha_1 > \frac{s_1(2m + m\theta + 1)}{m}$, the equilibrium point is a locally asymptotically stable node.

(ii) When $0 < s_2 < \frac{\alpha_2 s_1 - \alpha_1 \beta_2(\alpha_1 - \theta s_1 - s_1)(s_1 - m\alpha_1 + m\theta s_1)}{\alpha_1 - \theta s_1}$ and $\alpha_1 > \frac{s_1(2m + m\theta + 1)}{m}$, the equilibrium point is a locally asymptotically stable node.

(iii) When $s_2 = \frac{\alpha_2 s_1 - \alpha_1 \beta_2(\alpha_1 - \theta s_1 - s_1)(s_1 - m\alpha_1 + m\theta s_1)}{\alpha_1 - \theta s_1}$ and $\alpha_1 > \frac{s_1(2m + m\theta + 1)}{m}$, the equilibrium point is a saddle node of codimension 1.

Proof. The Jacobian matrix of system (7) at the equilibrium point

$$E_4 = (x_4, y_1^{*4}, y_2^{*4}) = \left(\frac{s_1}{\alpha_1 - \theta s_1}, \frac{\alpha_1(\alpha_1 - \theta s_1 - s_1)(s_1 - m\alpha_1 + m\theta s_1)}{(\alpha_1 - \theta s_1)^2}, 0\right)$$

is

$$J(E_4) = \begin{pmatrix} J_{11} & J_{12} & J_{13} \\ J_{21} & 0 & J_{23} \\ 0 & 0 & J_{33} \end{pmatrix},$$

where $J_{11} = -\frac{3s_1^2}{(\alpha_1 - \theta s_1)^2} + \frac{2(m+1)s_1}{\alpha_1 - \theta s_1} - m - \frac{(\alpha_1 - \theta s_1 - s_1)(s_1 - m\alpha_1 + m\theta s_1)(\alpha_1 - \theta s_1)}{\alpha_1}$, $J_{12} = -\frac{s_1}{\alpha_1}$,

$J_{13} = -\frac{s_1}{\alpha_1 - \theta s_1}$, $J_{21} = (\alpha_1 - \theta s_1 - s_1)(s_1 - m\alpha_1 + m\theta s_1)(\alpha_1 - \theta s_1)$,

$J_{23} = -\frac{\alpha_1 \beta_1(\alpha_1 - \theta s_1 - s_1)(s_1 - m\alpha_1 + m\theta s_1)}{\alpha_1 - \theta s_1}$,

$J_{33} = -s_2 - \frac{\alpha_1 \beta_2(\alpha_1 - \theta s_1 - s_1)(s_1 - m\alpha_1 + m\theta s_1) - \alpha_2 s_1}{\alpha_1 - \theta s_1}$,

and the characteristic equation is

$$|\lambda I - J(E_4)| = (\lambda - J_{33})(\lambda^2 - J_{11}\lambda - J_{12}J_{21}).$$

The eigenvalues of $J(E_4)$ is $\lambda_1 = J_{33}$, $\lambda_2$, and $\lambda_3$ are roots of the following equations:

$$\lambda^2 + (\frac{3s_1^2}{(\alpha_1 - \theta s_1)^2} - \frac{2(m+1)s_1}{\alpha_1 - \theta s_1} + m + \frac{(\alpha_1 - \theta s_1 - s_1)(s_1 - m\alpha_1 + m\theta s_1)(\alpha_1 - \theta s_1)}{\alpha_1})\lambda$$

$$+ \frac{s_1}{\alpha_1}(\alpha_1 - \theta s_1 - s_1)(s_1 - m\alpha_1 + m\theta s_1)(\alpha_1 - \theta s_1) = 0.$$

Since $y_1^{*4} = \frac{\alpha_1(\alpha_1 - \theta s_1 - s_1)(s_1 - m\alpha_1 + m\theta s_1)}{(\alpha_1 - \theta s_1)^2} > 0$, then

$$\lambda_2 \cdot \lambda_3 = \frac{s_1}{\alpha_1}(\alpha_1 - \theta s_1 - s_1)(s_1 - m\alpha_1 + m\theta s_1)(\alpha_1 - \theta s_1) > 0.$$

(i) Since $s_2 > \frac{\alpha_2 s_1 - \alpha_1 \beta_2(\alpha_1 - \theta s_1 - s_1)(s_1 - m\alpha_1 + m\theta s_1)}{\alpha_1 - \theta s_1}$, then

$$\lambda_1 = -s_2 - \frac{\alpha_1 \beta_2(\alpha_1 - \theta s_1 - s_1)(s_1 - m\alpha_1 + m\theta s_1) - \alpha_2 s_1}{\alpha_1 - \theta s_1} < 0,$$

Furthermore, since $x_4 = \frac{s_1}{\alpha_1 - \theta s_1} > 0$ and $\alpha_1 > \frac{s_1(2m + m\theta + 1)}{m}$, then

$$\lambda_2 + \lambda_3 = -\frac{3s_1^2}{(\alpha_1 - \theta s_1)^2} + \frac{2(m+1)s_1}{\alpha_1 - \theta s_1} - m - \frac{(\alpha_1 - \theta s_1 - s_1)(s_1 - m\alpha_1 + m\theta s_1)(\alpha_1 - \theta s_1)}{\alpha_1} < 0$$

hence the eigenvalues of $J(E_4)$ are $\lambda_1 < 0$, $\lambda_2 < 0$, and $\lambda_3 < 0$.

When $s_2 > \frac{\alpha_2 s_1 - \alpha_1 \beta_2(\alpha_1 - \theta s_1 - s_1)(s_1 - m\alpha_1 + m\theta s_1)}{\alpha_1 - \theta s_1}$ and $\alpha_1 > \frac{s_1(2m + m\theta + 1)}{m}$, the equilibrium point is locally asymptotically stable.

(ii) Since $0 < s_2 < \frac{\alpha_2 s_1 - \alpha_1 \beta_2(\alpha_1 - \theta s_1 - s_1)(s_1 - m\alpha_1 + m\theta s_1)}{\alpha_1 - \theta s_1}$, then

$$\lambda_1 = -s_2 - \frac{\alpha_1 \beta_2(\alpha_1 - \theta s_1 - s_1)(s_1 - m\alpha_1 + m\theta s_1) - \alpha_2 s_1}{\alpha_1 - \theta s_1} > 0,$$

and since $\alpha_1 > \frac{s_1(2m + m\theta + 1)}{m}$, then

$$\lambda_2 + \lambda_3 = -\frac{3s_1^2}{(\alpha_1 - \theta s_1)^2} + \frac{2(m+1)s_1}{\alpha_1 - \theta s_1} - m - \frac{(\alpha_1 - \theta s_1 - s_1)(s_1 - m\alpha_1 + m\theta s_1)(\alpha_1 - \theta s_1)}{\alpha_1} < 0.$$

Thus, the eigenvalues of $J(E_4)$ are $\lambda_1 > 0$, $\lambda_2 < 0$, and $\lambda_3 < 0$. Therefore, when

$$0 < s_2 < \frac{\alpha_2 s_1 - \alpha_1 \beta_2 (\alpha_1 - \theta s_1 - s_1)(s_1 - m\alpha_1 + m\theta s_1)}{\alpha_1 - \theta s_1} \text{ and } s_2 < \frac{(m+1)\alpha_2}{4}, \text{ the equilibrium point } E_4 \text{ is a}$$

saddle point.

(iii) When $s_2 = \dfrac{\alpha_2 s_1 - \alpha_1 \beta_2 (\alpha_1 - \theta s_1 - s_1)(s_1 - m\alpha_1 + m\theta s_1)}{\alpha_1 - \theta s_1}$ and

$$\alpha_2 > \frac{\alpha_1 \beta_2 (\alpha_1 - \theta s_1 - s_1)(s_1 - m\alpha_1 + m\theta s_1)}{s_1} > 0, \text{ we have } \lambda_1 = 0.$$

When $\alpha_1 > \dfrac{s_1(2m + m\theta + 1)}{m}$, we have

$$\lambda_2 + \lambda_3 = -\frac{3s_1^2}{(\alpha_1 - \theta s_1)^2} + \frac{2(m+1)s_1}{\alpha_1 - \theta s_1} - m - \frac{(\alpha_1 - \theta s_1 - s_1)(s_1 - m\alpha_1 + m\theta s_1)(\alpha_1 - \theta s_1)}{\alpha_1} < 0.$$

Thus, the eigenvalues of the matrix $J(E_4)$ are $\lambda_1 = 0$, $\lambda_2 < 0$, and $\lambda_3 < 0$.

Perform the transformation

$$\begin{cases} u = x - x_4 \\ v = y_1 - y_1^{*4}, \\ z = y_2 \end{cases}$$

by translating the equilibrium point to the origin and performing a Taylor expansion at the origin, we obtain:

$$\begin{cases} \dot{u} = a_{100}u + a_{010}v + a_{001}z + a_{200}u^2 + a_{110}uv + a_{101}uz + a_{300}u^3 + a_{210}u^2v + o(\|(u,v,z)\|^4) \\ \dot{v} = b_{100}u + b_{010}v + b_{001}z + b_{110}uv + b_{011}vz + b_{210}u^2v + o(\|(u,v,z)\|^4) \\ \dot{z} = c_{001}z + c_{101}uz + c_{011}vz + o(\|(u,v,z)\|^4) \end{cases}, \qquad (14)$$

where $a_{100} = -3x_4^2 + 2(m+1)x_4 - m - \dfrac{y_1^{*4}}{(1+\theta x_4)^2}$, $a_{010} = -\dfrac{x_4}{1+\theta x_4}$, $a_{001} = -x_4$,

$a_{200} = m + 1 - 3x_4 + \dfrac{y_1^{*4}}{(1+\theta x_4)^2}$, $a_{110} = -\dfrac{1}{(1+\theta x_4)^2}$, $a_{101} = -1$, $a_{300} = -1 - \dfrac{y_1^{*4}}{(1+\theta x_4)^4}$,

$a_{210} = \dfrac{1}{(1+\theta x_4)^3}$ , $b_{001} = \dfrac{\alpha_1 y_1^{*4}}{(1+\theta x_4)^2}$ , $b_{010} = -s_1 + \dfrac{\alpha_1 x_4}{1+\theta x_4}$ , $b_{001} = -\beta_1 y_1^{*4}$ , $b_{110} = \dfrac{\alpha_1}{(1+\theta x_4)^2}$ , $b_{011} = -\beta_1$ ,

$b_{210} = -\dfrac{\alpha_1}{(1+\theta x_4)^3}$ , $c_{001} = -s_2 + \alpha_2 x_4 - \beta_2 y_1^{*4}$ , $c_{101} = \alpha_2$ , $c_{011} = -\beta_2$ .

To transform the linearized matrix of system (14) into standard form, the following transformation is applied to system (14)

$$\begin{pmatrix} u \\ v \\ z \end{pmatrix} = \begin{pmatrix} w_{11} & w_{12} & w_{13} \\ w_{21} & w_{22} & w_{23} \\ w_{31} & w_{32} & w_{33} \end{pmatrix} \begin{pmatrix} u_1 \\ v_1 \\ z_1 \end{pmatrix},$$

where $w_{11} = 1$ , $w_{12} = 1$ , $w_{13} = 1$ , $w_{21} = \dfrac{\alpha_1(m-1)}{s_1}$ , $w_{22} = 0$ , $w_{23} = 0$ , $w_{31} = 0$ , $w_{32} = 0$ ,

$w_{33} = s_2 + m + 1 - \alpha_2$ .

Result:

$$\begin{cases} \dot{u}_1 = d_{200} u_1^2 + d_{110} u_1 v_1 + d_{101} u_1 z_1 + o(\|(u_1, v_1, z_1)\|^3) \\ \dot{v}_1 = \lambda_2 v_1 + e_{200} u_1^2 + e_{110} u_1 v_1 + e_{101} u_1 z_1 + e_{020} v_1^2 + e_{011} v_1 z_1 + e_{002} z_1^2 + o(\|(u_1, v_1, z_1)\|^3) , \quad (15) \\ \dot{z}_1 = \lambda_3 z_1 + f_{101} u_1 z_1 + f_{011} v_1 z_1 + f_{002} z_1^2 + o(\|(u_1, v_1, z_1)\|^3) \end{cases}$$

including $d_{200} = b_{110}$ , $d_{110} = b_{110}$ , $d_{101} = b_{110} + b_{011} w_{33}$ , $e_{200} = a_{200} + a_{110} w_{21} - b_{110}$ ,

$e_{110} = 2a_{200} + a_{110} w_{21} - b_{110}$, $e_{101} = 2a_{200} + a_{110} w_{21} + a_{101} w_{33} - b_{110} - b_{011} w_{33} - c_{101} - c_{011} w_{21}$, $e_{020} = a_{200}$,

$e_{011} = 2a_{200} + a_{101} w_{33} - c_{101}$ , $e_{002} = a_{200} + a_{101} w_{33} - c_{101}$ , $f_{101} = c_{101} + c_{011} w_{21}$ , $f_{011} = c_{101}$ , $f_{002} = c_{101}$.

Let the central manifold of (15) be $\begin{pmatrix} v_1 \\ z_1 \end{pmatrix} = \begin{pmatrix} h_1(u_1) \\ h_2(u_1) \end{pmatrix} = \begin{pmatrix} \tilde{\sigma}_1 u_1^2 + o(|u_1|^3) \\ \tilde{\sigma}_2 u_1^2 + o(|u_1|^3) \end{pmatrix}$ , then

$$\begin{pmatrix} 2\tilde{\sigma}_1 u_1 + \cdots \\ 2\tilde{\sigma}_2 u_1 + \cdots \end{pmatrix}(d_{200} u_1^2 + d_{110} u_1 \cdot \tilde{\sigma}_1 u_1^2 + d_{101} u_1 \cdot \tilde{\sigma}_2 u_1^2 + \cdots)$$
$$= \begin{pmatrix} \lambda_2 \cdot \tilde{\sigma}_1 u_1^2 + e_{200} u_1^2 + e_{110} u_1 \cdot \tilde{\sigma}_1 u_1^2 + e_{101} u_1 \cdot \tilde{\sigma}_2 u_1^2 + e_{020}(\tilde{\sigma}_1 u_1^2)^2 + e_{011} \cdot \tilde{\sigma}_1 u_1^2 \cdot \tilde{\sigma}_2 u_1^2 + e_{002}(\tilde{\sigma}_2 u_1^2)^2 + \cdots \\ \lambda_3 \cdot \tilde{\sigma}_2 u_1^2 + f_{101} u_1 \cdot \tilde{\sigma}_2 u_1^2 + f_{011} \cdot \tilde{\sigma}_1 u_1^2 \cdot \tilde{\sigma}_2 u_1^2 + f_{002}(\tilde{\sigma}_2 u_1^2)^2 + \cdots \end{pmatrix}.$$

By comparing the coefficients, we obtain $\lambda_2 \tilde{\sigma}_1 + e_{200} = 0$, $\lambda_3 \tilde{\sigma}_2 = 0$, so

$$\tilde{\sigma}_1 = -\dfrac{e_{200}}{\lambda_2}, \quad \tilde{\sigma}_2 = 0,$$

hence the central manifold is

$$\begin{pmatrix} v_1 \\ z_1 \end{pmatrix} = \begin{pmatrix} h_1(u_1) \\ h_2(u_1) \end{pmatrix} = \begin{pmatrix} -\dfrac{e_{200}}{\lambda_2} u_1^2 + o(|u_1|^3) \\ o(|u_1|^3) \end{pmatrix}.$$

System (15) is equivalent to the system near the origin:

$$\dot{u}_1 = g_{200} u_1^2 + g_{110} u_1^3 + o(\|(u_1, v_1, z_1)\|^3).$$

where $g_{200} = d_{200} = b_{110} = \dfrac{\alpha_1}{(1+\theta x_4)^2} > 0$, $g_{110} = -\dfrac{e_{200}}{\lambda_2} d_{110}$. By the theorem, we know that $E_3$ is a saddle-node point of codimension 1. □

**Theorem 7.** For the equilibrium point $E_5 = (x_5, y_1^{*5}, y_2^{*5})$, Let

$\beta_1 = \beta_2, \theta = 1, s_1 = s_2, m = \dfrac{2s_1}{\beta_1}$, and let

$$\alpha_1^* = \dfrac{5 s_1 \beta_1^2 \alpha_2 + 3 s_1^2 \beta_1 \alpha_2 + s_1^3 \alpha_2 + 5 \beta_1^3 \alpha_2 - 2 \beta_1^2 \alpha_2^2 - s_1 \beta_1 \alpha_2^2 - 3 s_1 \beta_1^3 - 3 s_1^2 \beta_1^2 - s_1^3 \beta_1 - 2 \beta_1^4}{\beta_1 (2 \beta_1 \alpha_2 - \beta_1^2 + s_1 \alpha_2)}, \text{ then}$$

there exists $\alpha_2^{*1}, \alpha_2^{*2}$ such that when $\alpha_2 \in (\alpha_2^{*1}, \alpha_2^{*2})$ and $\alpha_1 = \alpha_1^*$,

(i) If $\alpha_1 \neq \alpha_2 - \dfrac{(3 s_1 \beta_1^2 + 3 s_1^2 \beta_1 + s_1^3 + 2 \beta_1^3 - 2 \alpha_2 \beta_1^2 - \alpha_2 s_1 \beta_1)}{\beta_1^2 (\beta_1 + s_1)}$ and $\alpha_2 \neq \dfrac{3}{2} s_1$, the equilibrium point $E_5$ has one zero eigenvalue and two non-zero eigenvalues; this equilibrium point is a saddle-node point of codimension 1.

(ii) When $s_1 > 2\beta_1$ and $\alpha_2^* = \alpha_2(\beta_1, s_1)$, $s_1 \in (0, \delta)$, the equilibrium point $E_5$ has one zero eigenvalue and a pair of purely imaginary eigenvalues.

(iii) When $\tilde{\alpha}_2^* = \tilde{\alpha}_2(\beta_1, s_1)$, $s_1 \in (0, \tilde{\delta})$, the equilibrium point $E_5$ has two zero eigenvalues and one nonzero eigenvalue.

Proof. The Jacobian matrix of system (7) at the equilibrium point $E_5$ is

$$J(E_5) = \begin{pmatrix} \dfrac{s_1 (3 s_1 \beta_1^2 + 3 s_1^2 \beta_1 + s_1^3 + 2 \beta_1^3 - 2 \alpha_2 \beta_1^2 - \alpha_2 s_1 \beta_1)}{\beta_1^2 (\beta_1 + s_1)^2} & -\dfrac{s_1}{\beta_1 + s_1} & -\dfrac{s_1}{\beta_1} \\ \dfrac{\alpha_1 s_1 (\alpha_2 - \beta_1)}{(\beta_1 + s_1)^2} & \dfrac{s_1 (\alpha_1 - \alpha_2)}{\beta_1 + s_1} & -\dfrac{s_1 (\alpha_2 - \beta_1)}{\beta_1} \\ -\dfrac{\alpha_2 s_1 (\beta_1 + s_1 - \alpha_2)}{\beta_1 (\beta_1 + s_1)} & \dfrac{s_1 (\beta_1 + s_1 - \alpha_2)}{\beta_1 + s_1} & 0 \end{pmatrix}.$$

The characteristic equation is

$$|\lambda I - J(E_5)| = \lambda^3 - (\frac{s_1(3s_1\beta_1^2 + 3s_1^2\beta_1 + s_1^3 + 2\beta_1^3 - 2\alpha_2\beta_1^2 - \alpha_2 s_1\beta_1)}{\beta_1^2(\beta_1 + s_1)^2} + \frac{s_1(\alpha_1 - \alpha_2)}{\beta_1 + s_1})\lambda^2$$

$$+ (\frac{s_1^2(\alpha_1 - \alpha_2)(3s_1\beta_1^2 + 3s_1^2\beta_1 + s_1^3 + 2\beta_1^3 - 2\alpha_2\beta_1^2 - \alpha_2 s_1\beta_1)}{\beta_1^2(\beta_1 + s_1)^3} - \frac{\alpha_1 s_1^2(\alpha_2 - \beta_1)}{(\beta_1 + s_1)^3}$$

$$- \frac{\alpha_2 s_1^2(\beta_1 + s_1 - \alpha_2)}{\beta_1^2(\beta_1 + s_1)} + \frac{s_1^2(\alpha_2 - \beta_1)(\beta_1 + s_1 - \alpha_2)}{\beta_1(\beta_1 + s_1)})\lambda$$

$$- \frac{\alpha_2 s_1^3(\beta_1 + s_1 - \alpha_2)(\alpha_2 - \beta_1)}{\beta_1(\beta_1 + s_1)^3} + \frac{\alpha_2 s_1^3(\beta_1 + s_1 - \alpha_2)(\alpha_1 - \alpha_2)}{\beta_1^2(\beta_1 + s_1)^2}$$

$$- \frac{\alpha_1 s_1^3(\alpha_2 - \beta_1)(\beta_1 + s_1 - \alpha_2)}{\beta_1(\beta_1 + s_1)^3}$$

$$+ \frac{s_1^3(\beta_1 + s_1 - \alpha_2)(\alpha_2 - \beta_1)(3s_1\beta_1^2 + 3s_1^2\beta_1 + s_1^3 + 2\beta_1^3 - 2\alpha_2\beta_1^2 - \alpha_2 s_1\beta_1)}{\beta_1^2(\beta_1 + s_1)^3}.$$

(i) Since $x_5 = \frac{s_1}{\beta_1} > 0$, $y_1^{*5} = \frac{s_1(\alpha_2 - \beta_1)}{\beta_1^2} > 0$, and $y_2^{*5} = \frac{s_1(\alpha_2 - s_1 - \beta_1)}{\beta_1(\beta_1 + s_1)} > 0$,

then $\alpha_2 > s_1 + \beta_1$. So

$$2\beta_1\alpha_2 - \beta_1^2 + s_1\alpha_2 > 2\beta_1(s_1 + \beta_1) - \beta_1^2 + s_1(s_1 + \beta_1) = \beta_1^2 + s_1^2 + 3s_1\beta_1 > 0.$$

In addition, let

$$f(\alpha_2) = -(2\beta_1^2 + s_1\beta_1)\alpha_2^2 + (5s_1\beta_1^2 + 3s_1^2\beta_1 + s_1^3 + 5\beta_1^3)\alpha_2 - 3s_1\beta_1^3 - 3s_1^2\beta_1^2 - s_1^3\beta_1 - 2\beta_1^4,$$

there is

$$f(0) = -3s_1\beta_1^3 - 3s_1^2\beta_1^2 - s_1^3\beta_1 - 2\beta_1^4 < 0,$$

The axis of symmetry is $\alpha_2 = \frac{5s_1\beta_1^2 + 3s_1^2\beta_1 + s_1^3 + 5\beta_1^3}{2(2\beta_1^2 + s_1\beta_1)} > 0$.

Let

$$\Delta_{\alpha_2} = (5s_1\beta_1^2 + 3s_1^2\beta_1 + s_1^3 + 5\beta_1^3)^2 - 4(2\beta_1^2 + s_1\beta_1)(3s_1\beta_1^3 + 3s_1^2\beta_1^2 + s_1^3\beta_1 + 2\beta_1^4),$$

then $\exists\, 0 < \alpha_2^{*1} < \alpha_2^{*2}$ is a solution to $f(\alpha_2) = 0$,

where

$$\alpha_2^{*1} = \frac{(5s_1\beta_1^2 + 3s_1^2\beta_1 + s_1^3 + 5\beta_1^3) - \Delta_{\alpha_2}}{2(2\beta_1^2 + s_1\beta_1)} > 0, \quad \alpha_2^{*2} = \frac{(5s_1\beta_1^2 + 3s_1^2\beta_1 + s_1^3 + 5\beta_1^3) + \Delta_{\alpha_2}}{2(2\beta_1^2 + s_1\beta_1)} > 0.$$

When $\alpha_2 \in (\alpha_2^{*1}, \alpha_2^{*2})$, we have $f(\alpha_2) > 0$.

Thus, when $\alpha_2 \in (\alpha_2^{*1}, \alpha_2^{*2})$ and

$$\alpha_1 = \frac{5s_1\beta_1^2\alpha_2 + 3s_1^2\beta_1\alpha_2 + s_1^3\alpha_2 + 5\beta_1^3\alpha_2 - 2\beta_1^2\alpha_2^2 - s_1\beta_1\alpha_2^2 - 3s_1\beta_1^3 - 3s_1^2\beta_1^2 - s_1^3\beta_1 - 2\beta_1^4}{\beta_1(2\beta_1\alpha_2 - \beta_1^2 + s_1\alpha_2)} > 0,$$

we have

$$-\frac{\alpha_2 s_1^3(\beta_1 + s_1 - \alpha_2)(\alpha_2 - \beta_1)}{\beta_1(\beta_1 + s_1)^3} + \frac{\alpha_2 s_1^3(\beta_1 + s_1 - \alpha_2)(\alpha_1 - \alpha_2)}{\beta_1^2(\beta_1 + s_1)^2} - \frac{\alpha_1 s_1^3(\alpha_2 - \beta_1)(\beta_1 + s_1 - \alpha_2)}{\beta_1(\beta_1 + s_1)^3}$$

$$+\frac{s_1^3(\beta_1 + s_1 - \alpha_2)(\alpha_2 - \beta_1)(3s_1\beta_1^2 + 3s_1^2\beta_1 + s_1^3 + 2\beta_1^3 - 2\alpha_2\beta_1^2 - \alpha_2 s_1\beta_1)}{\beta_1^2(\beta_1 + s_1)^3} = 0.$$

If $\alpha_1 \neq \alpha_2 - \dfrac{(3s_1\beta_1^2 + 3s_1^2\beta_1 + s_1^3 + 2\beta_1^3 - 2\alpha_2\beta_1^2 - \alpha_2 s_1\beta_1)}{\beta_1^2(\beta_1 + s_1)}$, and

$$\frac{s_1^2(\alpha_1 - \alpha_2)(3s_1\beta_1^2 + 3s_1^2\beta_1 + s_1^3 + 2\beta_1^3 - 2\alpha_2\beta_1^2 - \alpha_2 s_1\beta_1)}{\beta_1^2(\beta_1 + s_1)^3} - \frac{\alpha_1 s_1^2(\alpha_2 - \beta_1)}{(\beta_1 + s_1)^3}$$

$$-\frac{\alpha_2 s_1^2(\beta_1 + s_1 - \alpha_2)}{\beta_1^2(\beta_1 + s_1)} + \frac{s_1^2(\alpha_2 - \beta_1)(\beta_1 + s_1 - \alpha_2)}{\beta_1(\beta_1 + s_1)} \neq 0.$$

In this case

$$|\lambda I - J(E_5)| = \lambda(\lambda^2 - (\frac{s_1(3s_1\beta_1^2 + 3s_1^2\beta_1 + s_1^3 + 2\beta_1^3 - 2\alpha_2\beta_1^2 - \alpha_2 s_1\beta_1)}{\beta_1^2(\beta_1 + s_1)^2} + \frac{s_1(\alpha_1 - \alpha_2)}{\beta_1 + s_1})\lambda$$

$$+(\frac{s_1^2(\alpha_1 - \alpha_2)(3s_1\beta_1^2 + 3s_1^2\beta_1 + s_1^3 + 2\beta_1^3 - 2\alpha_2\beta_1^2 - \alpha_2 s_1\beta_1)}{\beta_1^2(\beta_1 + s_1)^3} - \frac{\alpha_1 s_1^2(\alpha_2 - \beta_1)}{(\beta_1 + s_1)^3}$$

$$-\frac{\alpha_2 s_1^2(\beta_1 + s_1 - \alpha_2)}{\beta_1^2(\beta_1 + s_1)} + \frac{s_1^2(\alpha_2 - \beta_1)(\beta_1 + s_1 - \alpha_2)}{\beta_1(\beta_1 + s_1)})),$$

then the matrix $J(E_5)$ has an eigenvalue $\lambda_1 = 0$, $\lambda_2 \neq 0$, and $\lambda_3 \neq 0$.

Let

$$\Delta_\lambda = (\frac{s_1(3s_1\beta_1^2 + 3s_1^2\beta_1 + s_1^3 + 2\beta_1^3 - 2\alpha_2\beta_1^2 - \alpha_2 s_1\beta_1)}{\beta_1^2(\beta_1 + s_1)^2} + \frac{s_1(\alpha_1 - \alpha_2)}{\beta_1 + s_1})^2$$

$$-4(\frac{s_1^2(\alpha_1 - \alpha_2)(3s_1\beta_1^2 + 3s_1^2\beta_1 + s_1^3 + 2\beta_1^3 - 2\alpha_2\beta_1^2 - \alpha_2 s_1\beta_1)}{\beta_1^2(\beta_1 + s_1)^3}$$

$$-\frac{\alpha_1 s_1^2(\alpha_2 - \beta_1)}{(\beta_1 + s_1)^3} - \frac{\alpha_2 s_1^2(\beta_1 + s_1 - \alpha_2)}{\beta_1^2(\beta_1 + s_1)} + \frac{s_1^2(\alpha_2 - \beta_1)(\beta_1 + s_1 - \alpha_2)}{\beta_1(\beta_1 + s_1)}),$$

then

$$\lambda_2 = \frac{1}{2}\left( \frac{s_1(3s_1\beta_1^2 + 3s_1^2\beta_1 + s_1^3 + 2\beta_1^3 - 2\alpha_2\beta_1^2 - \alpha_2 s_1\beta_1)}{\beta_1^2(\beta_1 + s_1)^2} + \frac{s_1(\alpha_1 - \alpha_2)}{\beta_1 + s_1} - \sqrt{\Delta_\lambda} \right),$$

$$\lambda_3 = \frac{1}{2}\left( \frac{s_1(3s_1\beta_1^2 + 3s_1^2\beta_1 + s_1^3 + 2\beta_1^3 - 2\alpha_2\beta_1^2 - \alpha_2 s_1\beta_1)}{\beta_1^2(\beta_1 + s_1)^2} + \frac{s_1(\alpha_1 - \alpha_2)}{\beta_1 + s_1} + \sqrt{\Delta_\lambda} \right).$$

For the system

$$\begin{cases} \dfrac{dx}{dt} = x(1-x)\left(x - \dfrac{2s_1}{\beta_1}\right) - \dfrac{xy_1}{1+x} - xy_2 \\ \dfrac{dy_1}{dt} = -s_1 y_1 + \dfrac{\alpha_1 xy_1}{1+x} - \beta_1 y_1 y_2 \\ \dfrac{dy_2}{dt} = -s_2 y_2 + \alpha_2 xy_2 - \beta_2 y_1 y_2 \end{cases}.$$

Perform a transformation

$$\begin{cases} \tilde{u} = x - x_5 \\ \tilde{v} = y_1 - y_1^{*5}, \\ \tilde{z} = y_2 - y_2^{*5} \end{cases}$$

by applying a transformation to the system that shifts the equilibrium point to the origin and performing a Taylor expansion at the origin, we obtain:

$$\begin{cases} \dot{\tilde{u}} = a_{100}u + a_{010}v + a_{001}z + a_{200}u^2 + a_{110}uv + a_{101}uz + a_{300}u^3 + a_{210}u^2v + o(\|(u,v,z)\|^4) \\ \dot{\tilde{v}} = b_{100}u + b_{010}v + b_{001}z + b_{200}u^2 + b_{110}uv + b_{011}vz + b_{300}u^3 + b_{210}u^2v + o(\|(u,v,z)\|^4) \\ \dot{\tilde{z}} = c_{100}u + c_{010}v + c_{101}uz + c_{011}vz + o(\|(u,v,z)\|^4) \end{cases}, \quad (16)$$

where $a_{100} = -3x_5^2 + 2(m+1)x_5 - m - \dfrac{y_1^{*5}}{(1+x_5)^2} - y_2^{*5}$, $a_{010} = -\dfrac{x_5}{1+x_5}$, $a_{001} = -x_5$,

$a_{200} = m + 1 - 3x_5 - \dfrac{y_1^{*5}}{(1+x_5)^3}$, $a_{110} = -\dfrac{1}{(1+x_5)^2}$, $a_{101} = -1$, $a_{300} = -1 + \dfrac{y_1^{*5}}{(1+x_5)^4}$, $a_{210} = \dfrac{1}{(1+x_5)^3}$,

$b_{100} = \dfrac{\alpha_1 y_1^{*5}}{1+x_5}$, $b_{010} = -s_1 + \dfrac{\alpha_1 x_5}{1+x_5} - \beta_1 y_2^{*5}$, $b_{001} = \beta_1 y_1^{*5}$, $b_{200} = -\dfrac{\alpha_1 y_1^{*5}}{(1+x_5)^3}$, $b_{110} = \dfrac{\alpha_1}{(1+x_5)^2}$, $b_{011} = -\beta_1$,

$b_{300} = \dfrac{\alpha_1 y_1^{*5}}{3(1+x_5)^4}$, $b_{210} = -\dfrac{\alpha_1}{(1+x_5)^3}$, $c_{100} = \alpha_2 y_2^{*5}$, $c_{010} = -\beta_1 y_2^{*5}$, $c_{101} = \alpha_2$, $c_{011} = -\beta_1$.

To transform the linearized matrix of system (16) into standard form, the following transformation is applied to system (16)

$$\begin{pmatrix} u \\ v \\ z \end{pmatrix} = \begin{pmatrix} w_{11} & w_{12} & w_{13} \\ w_{21} & w_{22} & w_{23} \\ w_{31} & w_{32} & w_{33} \end{pmatrix} \begin{pmatrix} u_1 \\ v_1 \\ z_1 \end{pmatrix},$$

where $w_{11} = \dfrac{\alpha_2^2 (s_1 + \alpha_2 - \alpha_1)(\alpha_2 - s_1)^2}{s_1(\alpha_2^3 - 4\alpha_2^2 s_1 + 2\alpha_2 s_1^2 + 2\alpha_1 \alpha_2 s_1 - \alpha_1 s_1^2)}$, $w_{12} = 1$, $w_{13} = 1$

$w_{21} = -\dfrac{\alpha_2}{\alpha_2 - s_1} - \dfrac{\alpha_2 (2s_1 \alpha_2 - s_1^2)(s_1 + \alpha_2 - \alpha_1)}{\alpha_2^3 - 4\alpha_2^2 s_1 + 2\alpha_2 s_1^2 + 2\alpha_1 \alpha_2 s_1 - \alpha_1 s_1^2}$, $w_{22} = \dfrac{2s_1^3 \alpha_2^2 - s_1^4 \alpha_2 - \lambda_2 \alpha_2 (\alpha_2 - s_1)^2}{s_1}$,

$w_{23} = 0$,

$w_{31} = 1$, $w_{32} = \dfrac{2s_1^3 \alpha_2^2 - s_1^4 \alpha_2 - \lambda_2 \alpha_2 (\alpha_2 - s_1)^2}{s_1}$, $w_{33} = 0$.

We obtain

$$\begin{cases} \dot{u}_1 = d_{200} u_1^2 + d_{110} u_1 v_1 + d_{101} u_1 z_1 + o(\|(u_1, v_1, z_1)\|^3) \\ \dot{v}_1 = \lambda_2 v_1 + e_{200} u_1^2 + e_{110} u_1 v_1 + e_{101} u_1 z_1 + e_{020} v_1^2 + e_{011} v_1 z_1 + e_{002} z_1^2 + o(\|(u_1, v_1, z_1)\|^3) , \quad (17) \\ \dot{z}_1 = \lambda_3 z_1 + f_{200} u_1^2 + f_{110} u_1 v_1 + f_{101} u_1 z_1 + f_{020} v_1^2 + f_{011} v_1 z_1 + f_{002} z_1^2 + o(\|(u_1, v_1, z_1)\|^3) \end{cases}$$

where $d_{200} = c_{101} w_{11} + c_{011} w_{21}$, $d_{110} = c_{101} + c_{011} w_{22}$, $d_{101} = c_{101}$,

$e_{200} = \dfrac{1}{w_{22}}(b_{200} w_{11}^2 + b_{110} w_{11} w_{21} + b_{011} w_{21} + c_{101} w_{11}^2 w_{21} + c_{011} w_{21}^2 w_{11}$,
$- c_{101} w_{11}^2 w_{22} - c_{011} w_{11} w_{22} w_{21} - c_{101} w_{11} w_{21} - c_{011} w_{21}^2)$

$e_{110} = \dfrac{1}{w_{22}}(2b_{200} w_{11} + b_{110} w_{11} w_{22} + b_{110} w_{21} + b_{011} w_{22} + c_{101} w_{11} w_{21} + c_{011} w_{21} w_{11} w_{22}$,
$- c_{101} w_{11} w_{22} - c_{011} w_{11} w_{22}^2 - c_{101} w_{21} - c_{011} w_{21} w_{11} w_{22})$

$e_{101} = \dfrac{1}{w_{22}}(2b_{200} w_{11} + b_{110} w_{21} + c_{101} w_{11} w_{21} - c_{101} w_{11} w_{22} - c_{101} w_{21} - c_{101} w_{21} w_{11})$,

$e_{020} = \dfrac{1}{w_{22}}(b_{200} + b_{110} w_{22})$, $e_{011} = \dfrac{1}{w_{22}}(2b_{200} + b_{110} w_{22})$, $e_{002} = \dfrac{1}{w_{22}} b_{200}$,

$f_{200} = \dfrac{1}{w_{22}}(c_{101} w_{11} w_{21} + c_{011} w_{21}^2 - b_{200} w_{11}^2 - b_{110} w_{11} w_{21} - b_{011} w_{21} - c_{101} w_{11}^2 w_{21} - c_{011} w_{11} w_{21}^2)$,

$f_{110} = \dfrac{1}{w_{22}}(c_{101} w_{21} + c_{011} w_{21} w_{22} - 2b_{200} w_{11} - b_{110} w_{11} w_{22} - b_{110} w_{21} - b_{011} w_{22} - c_{101} w_{11} w_{21} - c_{011} w_{11} w_{21} w_{22})$,

$$f_{101} = 2a_{200}w_{11} + a_{110}w_{21} + a_{101} + \frac{1}{w_{22}}(c_{101}w_{21} - 2b_{200}w_{11} - b_{110}w_{21} - c_{101}w_{11}w_{21}),$$

$$f_{020} = a_{200} + a_{110}w_{22} - \frac{1}{w_{22}}(b_{200} + b_{110}w_{22}), f_{011} = 2a_{200} + a_{110}w_{22} - \frac{1}{w_{22}}(2b_{200} + b_{110}w_{22}),$$

$$f_{002} = a_{110}w_{23} - \frac{1}{w_{22}}b_{200}.$$

Let the central manifold of (17) be $\begin{pmatrix} v_1 \\ z_1 \end{pmatrix} = \begin{pmatrix} h_1(u_1) \\ h_2(u_1) \end{pmatrix} = \begin{pmatrix} \hat{\sigma}_1 u_1^2 + o(|u_1|^3) \\ \hat{\sigma}_2 u_1^2 + o(|u_1|^3) \end{pmatrix}$, then:

$$\begin{pmatrix} 2\hat{\sigma}_1 u_1 + \cdots \\ 2\hat{\sigma}_2 u_1 + \cdots \end{pmatrix}(d_{200}u_1^2 + d_{110}u_1 \cdot 2\hat{\sigma}_1 u_1 + d_{101}u_1 \cdot 2\hat{\sigma}_2 u_1 + \cdots)$$

$$= \begin{pmatrix} \lambda_2 \cdot 2\hat{\sigma}_1 u_1 + e_{200}u_1^2 + e_{110}u_1 \cdot 2\hat{\sigma}_1 u_1 + e_{101}u_1 \cdot 2\hat{\sigma}_2 u_1 \\ + e_{020}(2\hat{\sigma}_1 u_1)^2 + e_{011} \cdot 2\hat{\sigma}_1 u_1 \cdot 2\hat{\sigma}_2 u_1 + e_{002}(2\hat{\sigma}_2 u_1)^2 + \cdots \\ \lambda_3 \cdot 2\hat{\sigma}_2 u_1 + f_{200}u_1^2 + f_{110}u_1 \cdot 2\hat{\sigma}_1 u_1 + f_{101}u_1 \cdot 2\hat{\sigma}_2 u_1 \\ + f_{020}(2\hat{\sigma}_1 u_1)^2 + f_{011} \cdot 2\hat{\sigma}_1 u_1 \cdot z_1 + f_{002}(2\hat{\sigma}_2 u_1)^2 + \cdots \end{pmatrix}.$$

By comparing the coefficients, we obtain

$$\lambda_2 \hat{\sigma}_1 + e_{200} = 0, \quad \lambda_3 \hat{\sigma}_2 + f_{200} = 0, \quad \text{then } \hat{\sigma}_1 = -\frac{e_{200}}{\lambda_2}, \hat{\sigma}_2 = -\frac{f_{200}}{\lambda_3},$$

thus the central manifold is

$$\begin{pmatrix} v_1 \\ z_1 \end{pmatrix} = \begin{pmatrix} h_1(u_1) \\ h_2(u_1) \end{pmatrix} = \begin{pmatrix} -\frac{e_{200}}{\lambda_2}u_1^2 + o(|u_1|^3) \\ -\frac{f_{200}}{\lambda_3}u_1^2 + o(|u_1|^3) \end{pmatrix}.$$

System (16) is equivalent to the system near the origin:

$$\dot{u}_1 = g_{200}u_1^2 + g_{300}u_1^3 + o(\|u_1\|^4),$$

Where $g_{200} = d_{200}, g_{300} = -(\frac{d_{110}e_{200}}{\lambda_2} + \frac{d_{101}f_{200}}{\lambda_3})$.

For

$$d_{200} = c_{101}w_{11} + c_{011}w_{21}$$

$$= \alpha_2 + \frac{\alpha_2^{\,3}(s_1 + \alpha_2 - \alpha_1)(\alpha_2 - s_1)^2}{s_1(\alpha_2^{\,3} - 4\alpha_2^{\,2}s_1 + 2\alpha_2 s_1^{\,2} + 2\alpha_1\alpha_2 s_1 - \alpha_1 s_1^{\,2})} + \frac{\alpha_2(2s_1\alpha_2 - s_1^{\,2})(s_1 + \alpha_2 - \alpha_1)(\alpha_2 - s_1)}{\alpha_2^{\,3} - 4\alpha_2^{\,2}s_1 + 2\alpha_2 s_1^{\,2} + 2\alpha_1\alpha_2 s_1 - \alpha_1 s_1^{\,2}},$$

since $\alpha_2 > s_1$, we have $2s_1\alpha_2 - s_1^2 > 0$ and $\alpha_2 - s_1 > 0$. Furthermore, since $\alpha_2 \neq \frac{3}{2}s_1$,

we obtain $\alpha_2^3 - 4\alpha_2^2 s_1 + 2\alpha_2 s_1^2 + 2\alpha_1\alpha_2 s_1 - \alpha_1 s_1^2 \neq 0$.

Then $d_{200} = \alpha_2 > 0$, at this point, $g_{200} \neq 0$. By the theorem, $E_5$ is a saddle-node with codimension 1.

(ii) From the analysis in (i) above, it follows that when $\alpha_2 \in (\alpha_2^{*1}, \alpha_2^{*2})$ and

$\alpha_1 > 0$, we have

$$-\frac{\alpha_2 s_1^3(\beta_1 + s_1 - \alpha_2)(\alpha_2 - \beta_1)}{\beta_1(\beta_1 + s_1)^3} + \frac{\alpha_2 s_1^3(\beta_1 + s_1 - \alpha_2)(\alpha_1 - \alpha_2)}{\beta_1^2(\beta_1 + s_1)^2} - \frac{\alpha_1 s_1^3(\alpha_2 - \beta_1)(\beta_1 + s_1 - \alpha_2)}{\beta_1(\beta_1 + s_1)^3}$$
$$+ \frac{s_1^3(\beta_1 + s_1 - \alpha_2)(\alpha_2 - \beta_1)(3s_1\beta_1^2 + 3s_1^2\beta_1 + s_1^3 + 2\beta_1^3 - 2\alpha_2\beta_1^2 - \alpha_2 s_1\beta_1)}{\beta_1^2(\beta_1 + s_1)^3} = 0.$$

When $s_1 > 2\beta_1$, we have

$$f(\alpha_2)$$
$$= (-4s_1\beta_1^3 + s_1^2\beta_1^2 + s_1^3\beta_1 - 4\beta_1^4)\alpha_2^3$$
$$+ (16\beta_1^5 + 15s_1\beta_1^4 + 9s_1^2\beta_1^3 + 5s_1^3\beta_1^2 + s_1^4\beta_1 + 4\beta_1^4 + 4s_1\beta_1^3 + s_1^2\beta_1^2)\alpha_2^2$$
$$- (17\beta_1^6 + 29s_1\beta_1^5 + 38s_1^2\beta_1^4 + 38s_1^3\beta_1^3 + 17s_1^4\beta_1^2 + 6s_1^5\beta_1 + s_1^6$$
$$+ 9s_1\beta_1^4 + 9s_1^2\beta_1^3 + 5s_1^3\beta_1^2 + s_1^4\beta_1 + 6\beta_1^5)\alpha_2$$
$$+ 15s_1\beta_1^6 + 24s_1^2\beta_1^5 + 23s_1^3\beta_1^4 + 15s_1^4\beta_1^3 + 6s_1^5\beta_1^2 + s_1^6\beta_1 + 6\beta_1^7 + \beta_1^6 > 0.$$

Thus

$$(\frac{s_1^2(\alpha_1 - \alpha_2)(3s_1\beta_1^2 + 3s_1^2\beta_1 + s_1^3 + 2\beta_1^3 - 2\alpha_2\beta_1^2 - \alpha_2 s_1\beta_1)}{\beta_1^2(\beta_1 + s_1)^3} - \frac{\alpha_1 s_1^2(\alpha_2 - \beta_1)}{(\beta_1 + s_1)^3}$$
$$- \frac{\alpha_2 s_1^2(\beta_1 + s_1 - \alpha_2)}{\beta_1^2(\beta_1 + s_1)} + \frac{s_1^2(\alpha_2 - \beta_1)(\beta_1 + s_1 - \alpha_2)}{\beta_1(\beta_1 + s_1)}) > 0.$$

Let

$$g(\alpha_2, s_1, \beta_1)$$
$$= -(8\beta_1^3 + 10s_1\beta_1^2 + 3s_1^2\beta_1)\alpha_2^2$$
$$+ (12\beta_1^4 + 20s_1\beta_1^3 + 17s_1^2\beta_1^2 + 9s_1^3\beta_1 + 2s_1^4)\alpha_2$$
$$- (8s_1\beta_1^4 + 9s_1^2\beta_1^3 + 5s_1^3\beta_1^2 + s_1^4\beta_1 + 4\beta_1^5),$$

since $y_1^{*5} = \frac{s_1(\alpha_2 - \beta_1)}{\beta_1^2} > 0$, we have $\alpha_2 \neq \beta_1$, so

$$\left.\frac{\partial g}{\partial \alpha_2}\right|_{s_1=0} = -8\beta_1^3\alpha_2^2 + 12\beta_1^4\alpha_2 - 4\beta_1^5 \neq 0.$$

By the Implicit Function Theorem, there exists a unique function $\alpha_2^* = \alpha_2(\beta_1, s_1)$, $s_1 \in (0, \delta)$ such that $g(\alpha_2^*, s_1, \beta_1) = 0$.

Thus, when $\alpha_2 = \alpha_2^*$, we have

$$\frac{s_1(3s_1\beta_1^2 + 3s_1^2\beta_1 + s_1^3 + 2\beta_1^3 - 2\alpha_2\beta_1^2 - \alpha_2 s_1 \beta_1)}{\beta_1^2(\beta_1+s_1)^2} + \frac{s_1(\alpha_1 - \alpha_2)}{\beta_1 + s_1} = 0,$$

The characteristic equation is

$$|\lambda I - J(E_5)| = \lambda(\lambda^2 + (\frac{s_1^2(\alpha_1 - \alpha_2)(3s_1\beta_1^2 + 3s_1^2\beta_1 + s_1^3 + 2\beta_1^3 - 2\alpha_2\beta_1^2 - \alpha_2 s_1 \beta_1)}{\beta_1^2(\beta_1+s_1)^3}$$
$$- \frac{\alpha_1 s_1^2(\alpha_2 - \beta_1)}{(\beta_1+s_1)^3} - \frac{\alpha_2 s_1^2(\beta_1+s_1-\alpha_2)}{\beta_1^2(\beta_1+s_1)} + \frac{s_1^2(\alpha_2-\beta_1)(\beta_1+s_1-\alpha_2)}{\beta_1(\beta_1+s_1)})).$$

In this case, there is one zero eigenvalue and two purely imaginary eigenvalues.

(iii) From (i), we know that when $\alpha_2 \in (\alpha_2^{*1}, \alpha_2^{*2})$ and $\alpha_1 > 0$, we have

$$-\frac{\alpha_2 s_1^3(\beta_1+s_1-\alpha_2)(\alpha_2-\beta_1)}{\beta_1(\beta_1+s_1)^3} + \frac{\alpha_2 s_1^3(\beta_1+s_1-\alpha_2)(\alpha_1-\alpha_2)}{\beta_1^2(\beta_1+s_1)^2} - \frac{\alpha_1 s_1^3(\alpha_2-\beta_1)(\beta_1+s_1-\alpha_2)}{\beta_1(\beta_1+s_1)^3}$$
$$+ \frac{s_1^3(\beta_1+s_1-\alpha_2)(\alpha_2-\beta_1)(3s_1\beta_1^2 + 3s_1^2\beta_1 + s_1^3 + 2\beta_1^3 - 2\alpha_2\beta_1^2 - \alpha_2 s_1 \beta_1)}{\beta_1^2(\beta_1+s_1)^3} = 0.$$

Let

$$g(\alpha_2, s_1, \beta_1)$$
$$= (-4s_1\beta_1^3 + s_1^2\beta_1^2 + s_1^3\beta_1 - 4\beta_1^4)\alpha_2^3$$
$$+ (16\beta_1^5 + 15s_1\beta_1^4 + 9s_1^2\beta_1^3 + 5s_1^3\beta_1^2 + s_1^4\beta_1 + 4\beta_1^4 + 4s_1\beta_1^3 + s_1^2\beta_1^2)\alpha_2^2$$
$$- (17\beta_1^6 + 29s_1\beta_1^5 + 38s_1^2\beta_1^4 + 38s_1^3\beta_1^3 + 17s_1^4\beta_1^2 + 6s_1^5\beta_1 + s_1^6$$
$$+ 9s_1\beta_1^4 + 9s_1^2\beta_1^3 + 5s_1^3\beta_1^2 + s_1^4\beta_1 + 6\beta_1^5)\alpha_2$$
$$+ 15s_1\beta_1^6 + 24s_1^2\beta_1^5 + 23s_1^3\beta_1^4 + 15s_1^4\beta_1^3 + 6s_1^5\beta_1^2 + s_1^6\beta_1 + 6\beta_1^7 + \beta_1^6.$$

Furthermore,

$$\left.\frac{\partial g}{\partial \alpha_2}\right|_{\alpha_2=\frac{8+\sqrt{13}}{6}\beta_1, s_1=0}$$

$$= -12\beta_1^4\alpha_2^2 + 8(4\beta_1^5 + \beta_1^4)\alpha_2 - (17\beta_1^6 + 6\beta_1^5)\Big|_{\alpha_2=\frac{8+\sqrt{13}}{6}\beta_1, s_1=0}$$

$$= \frac{14+4\sqrt{13}}{3}\beta_1^5 \neq 0.$$

By the Implicit Function Theorem, there exists a unique function $\tilde{\alpha}_2^* = \tilde{\alpha}_2(\beta_1, s_1)$, $s_1 \in (0, \tilde{\delta})$ such that $g(\tilde{\alpha}_2^*, s_1, \beta_1) = 0$.

Thus, when $\alpha_2 = \tilde{\alpha}_2^*$, we have

$$\frac{s_1^2(\alpha_1 - \alpha_2)(3s_1\beta_1^2 + 3s_1^2\beta_1 + s_1^3 + 2\beta_1^3 - 2\alpha_2\beta_1^2 - \alpha_2 s_1\beta_1)}{\beta_1^2(\beta_1 + s_1)^3} - \frac{\alpha_1 s_1^2(\alpha_2 - \beta_1)}{(\beta_1 + s_1)^3}$$

$$- \frac{\alpha_2 s_1^2(\beta_1 + s_1 - \alpha_2)}{\beta_1^2(\beta_1 + s_1)} + \frac{s_1^2(\alpha_2 - \beta_1)(\beta_1 + s_1 - \alpha_2)}{\beta_1(\beta_1 + s_1)} = 0.$$

If $\alpha_1 \neq \alpha_2 - \dfrac{(3s_1\beta_1^2 + 3s_1^2\beta_1 + s_1^3 + 2\beta_1^3 - 2\alpha_2\beta_1^2 - \alpha_2 s_1\beta_1)}{\beta_1^2(\beta_1 + s_1)}$, then

$$|\lambda I - J(E_5)| = \lambda^2 \left(\lambda - \left(\frac{s_1(3s_1\beta_1^2 + 3s_1^2\beta_1 + s_1^3 + 2\beta_1^3 - 2\alpha_2\beta_1^2 - \alpha_2 s_1\beta_1)}{\beta_1^2(\beta_1 + s_1)^2} + \frac{s_1(\alpha_1 - \alpha_2)}{\beta_1 + s_1}\right)\right),$$

then the matrix $J(E_5)$ has eigenvalues $\lambda_1 = 0$, $\lambda_2 = 0$,

$$\lambda_3 = \frac{s_1(3s_1\beta_1^2 + 3s_1^2\beta_1 + s_1^3 + 2\beta_1^3 - 2\alpha_2\beta_1^2 - \alpha_2 s_1\beta_1)}{\beta_1^2(\beta_1 + s_1)^2} + \frac{s_1(\alpha_1 - \alpha_2)}{\beta_1 + s_1} \neq 0.$$

By the central manifold theorem, let the central manifold of (16) be

$$z = \sigma_{10}u + \sigma_{01}v + \sigma_{20}u^2 + \sigma_{11}uv + \sigma_{02}v^2 + (\|(u,v)\|^3). \tag{18}$$

Taking the derivative of the central manifold with respect to $t$ at both ends simultaneously yields

$$\frac{dz}{dt} = \sigma_{10}\frac{du}{dt} + \sigma_{01}\frac{dv}{dt} + 2\sigma_{20}u\frac{du}{dt} + \sigma_{11}u\frac{dv}{dt} + \sigma_{11}v\frac{du}{dt} + 2\sigma_{02}v\frac{dv}{dt} + \cdots$$

$$=\sigma_{10}(a_{100}u + a_{010}v + a_{001}(\sigma_{10}u + \sigma_{01}v + \sigma_{20}u^2 + \sigma_{11}uv + \sigma_{02}v^2 + \cdots) + a_{200}u^2 + a_{110}uv$$
$$+ a_{101}u(\sigma_{10}u + \sigma_{01}v + \sigma_{20}u^2 + \sigma_{11}uv + \sigma_{02}v^2 + \cdots) + \cdots)$$
$$+\sigma_{01}(b_{100}u + b_{010}v + b_{001}(\sigma_{10}u + \sigma_{01}v + \sigma_{20}u^2 + \sigma_{11}uv + \sigma_{02}v^2 + \cdots) + b_{200}u^2 + b_{110}uv$$
$$+ b_{011}v(\sigma_{10}u + \sigma_{01}v + \sigma_{20}u^2 + \sigma_{11}uv + \sigma_{02}v^2 + \cdots) + \cdots)$$
$$+ 2\sigma_{20}u(a_{100}u + a_{010}v + a_{001}(\sigma_{10}u + \sigma_{01}v + \sigma_{20}u^2 + \sigma_{11}uv + \sigma_{02}v^2 + \cdots) + a_{200}u^2 + a_{110}uv$$
$$+ a_{101}u(\sigma_{10}u + \sigma_{01}v + \sigma_{20}u^2 + \sigma_{11}uv + \sigma_{02}v^2 + \cdots) + \cdots)$$
$$+ \sigma_{11}u(b_{100}u + b_{010}v + b_{001}(\sigma_{10}u + \sigma_{01}v + \sigma_{20}u^2 + \sigma_{11}uv + \sigma_{02}v^2 + \cdots) + b_{200}u^2 + b_{110}uv$$
$$+ b_{011}v(\sigma_{10}u + \sigma_{01}v + \sigma_{20}u^2 + \sigma_{11}uv + \sigma_{02}v^2 + \cdots) + \cdots)$$
$$+ \sigma_{11}v(a_{100}u + a_{010}v + a_{001}(\sigma_{10}u + \sigma_{01}v + \sigma_{20}u^2 + \sigma_{11}uv + \sigma_{02}v^2 + \cdots) + a_{200}u^2 + a_{110}uv$$
$$+ a_{101}u(\sigma_{10}u + \sigma_{01}v + \sigma_{20}u^2 + \sigma_{11}uv + \sigma_{02}v^2 + \cdots) + \cdots)$$
$$+ 2\sigma_{02}v(b_{100}u + b_{010}v + b_{001}(\sigma_{10}u + \sigma_{01}v + \sigma_{20}u^2 + \sigma_{11}uv + \sigma_{02}v^2 + \cdots) + b_{200}u^2 + b_{110}uv$$
$$+ b_{011}v(\sigma_{10}u + \sigma_{01}v + \sigma_{20}u^2 + \sigma_{11}uv + \sigma_{02}v^2 + \cdots) + \cdots)$$
$$+ \cdots.$$

Furthermore, substituting the central manifold (17) into the differential equation for the variable $z$ with respect to $t$ yields:

$$\frac{dz}{dt} = c_{100}u + c_{010}v + c_{101}u(\sigma_{10}u + \sigma_{01}v + \sigma_{20}u^2 + \sigma_{11}uv + \sigma_{02}v^2 + \cdots)$$
$$+ c_{011}v(\sigma_{10}u + \sigma_{01}v + \sigma_{20}u^2 + \sigma_{11}uv + \sigma_{02}v^2 + \cdots)$$
$$+ \cdots.$$

By comparing the coefficients, we obtain:

$$a_{100}\sigma_{10} + a_{001}\sigma_{10}^2 + b_{100}\sigma_{01} + b_{001}\sigma_{10}\sigma_{01} = c_{100},$$

$$a_{010}\sigma_{10} + a_{100}\sigma_{10}\sigma_{01} + b_{010}\sigma_{01} + b_{001}\sigma_{01}^2 = c_{010},$$

$$a_{001}\sigma_{20} + a_{200}\sigma_{10} + a_{101}\sigma_{10}^2 + b_{001}\sigma_{01} + b_{200}\sigma_{01} + 2a_{100}\sigma_{20} + 2a_{001}\sigma_{20}\sigma_{10} + b_{100}\sigma_{11} + b_{001}\sigma_{10}\sigma_{11} = c_{101}\sigma_{10},$$

$$a_{001}\sigma_{11} + a_{110}\sigma_{10} + a_{101}\sigma_{10}\sigma_{01} + b_{001}\sigma_{01}\sigma_{11} + b_{110}\sigma_{01} + b_{011}\sigma_{10}\sigma_{01} + 2a_{010}\sigma_{20} + 2a_{001}\sigma_{01}\sigma_{20}$$
$$+ b_{010}\sigma_{11} + b_{001}\sigma_{01}\sigma_{11} + a_{100}\sigma_{11} + a_{001}\sigma_{10}\sigma_{11} + 2b_{100}\sigma_{02} + 2b_{001}\sigma_{10}\sigma_{02} = c_{101}\sigma_{01} + c_{011}\sigma_{10},$$

$$a_{001}\sigma_{02} + b_{001}\sigma_{01}\sigma_{02} + b_{011}\sigma_{01}^2 + a_{010}\sigma_{11} + a_{001}\sigma_{01}\sigma_{11} + 2b_{010}\sigma_{20} + 2b_{001}\sigma_{01}\sigma_{02} = c_{011}\sigma_{01}.$$

Let the solution to the above system of equations be $(\sigma_{10}^*, \sigma_{01}^*, \sigma_{20}^*, \sigma_{11}^*, \sigma_{02}^*)$. Numerical simulations show that the system has a unique solution. Then the central manifold is

$$z = \sigma_{10}^* u + \sigma_{01}^* v + \sigma_{20}^* u^2 + \sigma_{11}^* uv + \sigma_{02}^* v^2 + (\|(u,v)\|^3).$$

Substituting the obtained central manifold back into system (15), it is equivalent to the following two-dimensional system near the origin:

$$\begin{cases} \dot{u}_1 = d_{10}u_1 + d_{01}v_1 + d_{20}u_1^2 + d_{11}u_1v_1 + d_{20}v_1^2 + d_{30}u_1^3 + d_{21}u_1^2v_1 + d_{12}u_1v_1^2 + o(\|(u_1,v_1)\|^4) \\ \dot{v}_1 = e_{10}u_1 + e_{01}v_1 + e_{20}u_1^2 + e_{11}u_1v_1 + e_{20}v_1^2 + e_{30}u_1^3 + e_{21}u_1^2v_1 + e_{12}u_1v_1^2 + e_{03}v_1^3 + o(\|(u_1,v_1)\|^4) \end{cases}, \text{where}$$

$$d_{10} = a_{100} + a_{001}\sigma_{10}^*, \; d_{01} = a_{010} + a_{001}\sigma_{01}^*, \; d_{20} = a_{001}\sigma_{20}^* u^2 + a_{200} + a_{101}\sigma_{10}^*,$$

$$d_{11} = a_{001}\sigma_{11}^* + a_{110} + a_{101}\sigma_{01}^*, \; d_{02} = a_{001}\sigma_{02}^*, \; d_{30} = a_{101}\sigma_{20}^* + a_{300}, \; d_{21} = a_{101}\sigma_{11}^* + a_{210},$$

$$d_{12} = a_{101}\sigma_{02}^*, \; e_{10} = b_{100}u + b_{001}\sigma_{10}^*, \; e_{01} = b_{010} + b_{001}\sigma_{01}^*, \; e_{20} = a_{001}\sigma_{20}^* + a_{200} + a_{101}\sigma_{10}^*$$

$$e_{11} = a_{001}\sigma_{11}^* + a_{110} + a_{101}\sigma_{01}^*, \; e_{02} = a_{001}\sigma_{02}^* v^2, \; e_{30} = b_{300}u^3, \; e_{21} = b_{011}\sigma_{20}^* + b_{210},$$

$$e_{12} = b_{011}\sigma_{11}^*, \; e_{03} = b_{011}\sigma_{02}^*. \qquad \square$$

## 3 Numerical Simulation

For the equilibrium point $E_1$, values satisfying the conditions of Theorem 3 were selected. For example, taking $m = 0.2$, $\alpha_1 = 0.2$, $\alpha_2 = 0.1$, $\beta_1 = 0.15$, $\theta = 0.5$, $s_1 = 0.2$, $s_2 = 0.3$ and $\beta_2 = 0.12$, the following graph shows the time-dependent changes in the population densities of the three species.

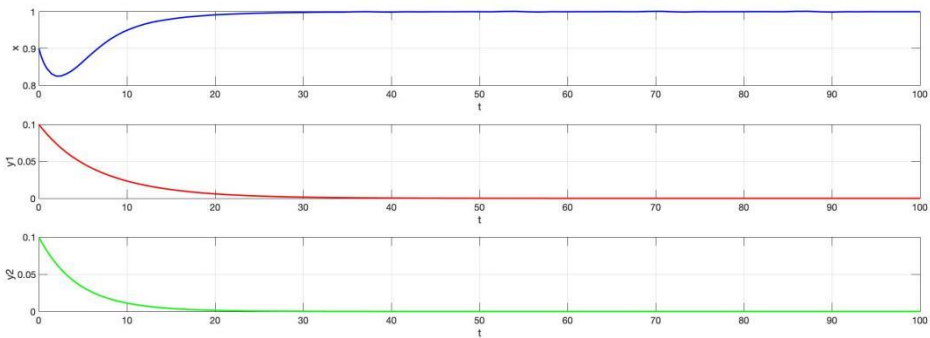

Figure 1: Population density curves of three species when $\dfrac{\alpha_1}{1+\theta} < s_1$ and $\alpha_2 < s_2$

From Figure 1, it can be seen that for the given initial conditions $(0.9, 0.1, 0.1)$, the population density of the prey population $x$ first decreases and then increases, eventually converging to a steady-state value; the population density of the predator population $y_1$ continues to decrease until it reaches zero; simultaneously, compared to the predator population $y_1$, the population density of the predator population $y_2$ decreases at a faster rate until it reaches zero. This trend confirms that when $\frac{\alpha_1}{1+\theta} < s_1$ and $\alpha_2 < s_2$, the equilibrium point $E_1$ is locally stable, consistent with the analysis.

For the equilibrium point $E_3$, select values that satisfy the conditions of Theorem 5. For example, taking $m = 0.02, \alpha_1 = 0.2, \alpha_2 = 0.1, \beta_1 = 0.15, \theta = 0.5$, $s_1 = 0.08, s_2 = 0.04$, and $\beta_2 = 0.12$ yields the following graph showing the relationship between the population densities of the three species over time.

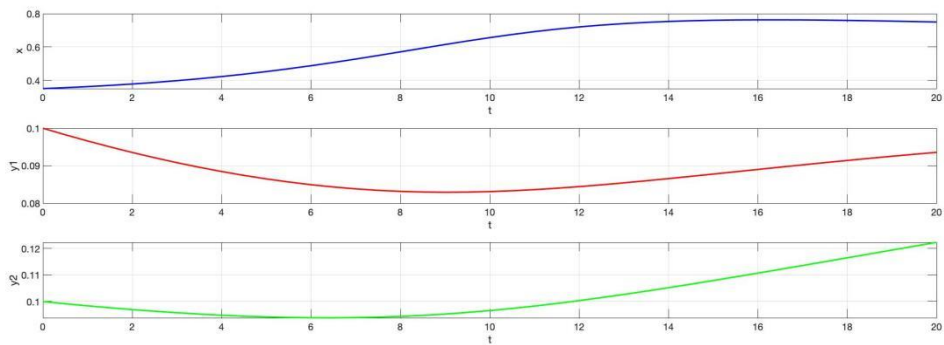

Figure 2: Population density curves of three species when condition (1) of Theorem 5 is satisfied

From Figure 2, it can be seen that for the given initial conditions $(0.35, 0.1, 0.1)$, the population density of the prey species $x$ will rise and eventually converge to a steady-state value; the population density of the predator species $y_1$ will first decrease and then increase before converging to a steady-state value; simultaneously, the population density of the predator species $y_2$ will also first decrease and then increase before converging to a steady-state value. This trend confirms that when Theorem 5 is satisfied, the equilibrium point $E_3$ is locally stable, consistent with the analysis.

## 4 Conclusions

The main contribution of this paper is to study a three-species predator-prey model based on the classical Lotka-Volterra model, considering the Allee effect in the prey population and interspecific competition between the two predator populations. We discuss the existence of various boundary equilibrium points in the model. Under sufficient conditions, determine the existence of internal equilibrium points: there may be zero, one, or two internal equilibrium points;Furthermore, the stability of boundary equilibrium points under different parameter conditions is analyzed; these may be stable or unstable nodes, saddles, saddle-nodes, or cusp points of dimension 2. For internal equilibrium points, their eigenvalue conditions are analyzed, and the parameter conditions required for cases with one zero eigenvalue and two non-zero eigenvalues, one zero eigenvalue and a pair of purely imaginary eigenvalues, and two zero eigenvalues and one non-zero eigenvalue are provided.When an internal equilibrium point has one zero eigenvalue and a pair of purely imaginary eigenvalues, the system may undergo a fold-Hopf bifurcation; when an internal equilibrium point has two zero eigenvalues and one non-zero eigenvalue, the system may undergo a B-T bifurcation. Given the computational complexity and the complexity of the notation, further research and analysis may be conducted in subsequent studies.